\documentclass{article}

    \PassOptionsToPackage{numbers, compress}{natbib}

 \usepackage[preprint]{neurips_2026}

\usepackage[utf8]{inputenc}
\usepackage[T1]{fontenc}
\usepackage{hyperref}
\usepackage{url}
\usepackage{booktabs}
\usepackage{amsfonts}
\usepackage{nicefrac}
\usepackage{microtype}
\usepackage{xcolor}
\usepackage{array}
\usepackage{amsmath}
\usepackage{amsthm}
\usepackage{mathtools}
\usepackage{graphicx}
\graphicspath{{figures/}}
\usepackage{subcaption}
\usepackage{algorithm}
\usepackage{algpseudocode}
\usepackage{bm}

\usepackage{algorithm}

\newcommand{\KL}{\operatorname{KL}}

\newcommand{\R}{\mathbb{R}}
\newcommand{\e}{\varepsilon}
\newcommand{\OT}{\text{OT}}
\newcommand{\EOT}{\text{EOT}}

\newcommand{\one}{\mathbf{1}}

\DeclareMathOperator*{\argmin}{arg\,min}

\theoremstyle{plain}
\newtheorem{theorem}{Theorem}
\newtheorem{proposition}[theorem]{Proposition}
\newtheorem{corollary}[theorem]{Corollary}
\newtheorem{lemma}[theorem]{Lemma}

\theoremstyle{definition}
\newtheorem{example}[theorem]{Example}

\theoremstyle{definition}
\newtheorem{remark}[theorem]{Remark}

\title{Convergence of potential-preserving annealed Sinkhorn algorithm}

\author{Dawid Lipinski, Maciej Buze\\
  MARS: Mathematics for AI in Real-world Systems\\
  School of Mathematical Sciences\\
  Lancaster University\\
  Lancaster, LA1 4YF, United Kingdom\\
  \texttt{\{d.a.lipinski, m.buze\}@lancaster.ac.uk}
}

\begin{document}

\maketitle

\begin{abstract}
The annealed Sinkhorn algorithm is an iterative method for approximating solutions of optimal transport problems, in which the entropic regularisation parameter is progressively decreased.
    As this parameter decreases, one may preserve either the multiplicative scaling factors or the corresponding dual potentials.
    We study the potential-preserving iterates under positive, non-increasing regularisation schedules tending to zero.
    We prove that, if the schedule is non-summable, the iterates admit a subsequence converging to a solution of the corresponding optimal transport problem.
    The non-summability condition is sharp, as for every summable schedule in this class we construct an example in which the iterates converge to a plan violating the prescribed source marginal.
    We further establish convergence of the entire sequence if the relative inverse-temperature increments are also square-summable, and show by example that non-summability alone does not guarantee full convergence.
    Our convergence conditions are satisfied by linear inverse-temperature schedules, under which recent work has shown that the scaling-preserving iterates fail to converge to a solution of the optimal transport problem.
    We also derive convergence rates for polynomial inverse-temperature schedules up to and including linear growth, extending the marginal KL bound obtained by mirror-descent analysis at fixed temperature.
    Our best rate for the transport cost gap after rounding improves on that established for scaling-preserving annealing. 
    Simple numerical illustrations accompany these convergence results.

\end{abstract}

\section{Introduction}
\label{sec:Intro}

Optimal transport (OT) compares probability distributions through their geometry and simultaneously provides a coupling between their samples. Its applications range from distribution matching and generative modelling \citep{lipman2023flow, tong_improving_2024, patrini_sinkhorn_2020}, domain adaptation \citep{courty_optimal_2017}, graph and shape correspondence \citep{xu_gromov-wasserstein_2019,vayer_optimal_2019} to genetics \citep{schiebinger_optimal-transport_2019,passemiers_dagip_2025} and biology \citep{bunne_learning_2023,klein_mapping_2025,bryan_accurate_2026}, as well as physics and weather forecasting
\citep{komiske_metric_2019,francis_examining_2025}. Entropic regularisation makes discrete OT computationally tractable by adding an entropy term with the \emph{temperature} prefactor $\varepsilon > 0$, thus reducing the problem to the matrix-scaling iterations of Sinkhorn's algorithm \citep{cuturi_sinkhorn_2013,altschuler_near-linear_2017}. Choosing $\e$, however, is challenging: strong regularisation introduces substantial bias, whereas weak regularisation leads to slow convergence \citep{feydy_interpolating_2019}.
Annealed Sinkhorn addresses this trade-off by gradually decreasing the temperature along $(\varepsilon_t)_{t\in \mathbb N}$ and warm-starting each regularisation level from the preceding iterate \citep{cuturi_sinkhorn_2013, schmitzer_stabilized_2019, chizat_annealed_2024}. This strategy has been explicitly adopted in recent large-scale OT and Gromov-Wasserstein solvers \citep{houry_gromov-wasserstein_2026,ye_flashsinkhorn_2026}. A key implementation choice is how to transfer the iterates between successive temperatures: one may preserve either the multiplicative scalings or the corresponding dual potentials, with the latter being more commonly used in practice \citep{flamary_pot_2021,feydy_interpolating_2019}.
Despite its practical use, including in the stabilised $\varepsilon$-scaling algorithm of \citep{schmitzer_stabilized_2019}, the convergence of potential-preserving annealing remains unestablished.
In contrast, the scaling-preserving annealed Sinkhorn algorithm has recently been analysed under concave schedules \citep{chizat_annealed_2024}: its iterates converge to a solution of the corresponding OT problem precisely when the inverse temperature $\beta_t=\varepsilon_t^{-1}\to\infty$ and $\beta_t-\beta_{t-1}\to0$, whereas under a linear schedule $\beta_t\sim t$ it converges instead to a solution of a semi-unbalanced problem \citep{liero_optimal_2018}. The analysis in \citep{chizat_annealed_2024} studies a computationally tractable proxy for scaling-preserving annealed Sinkhorn and shows that its best convergence rate is \(O(t^{-1/2})\). To reduce the bias introduced by scaling preservation, \citep{chizat_annealed_2024} also proposes a debiased variant of annealed Sinkhorn. Numerical experiments suggest that this modification permits faster annealing, although its convergence to a solution of the balanced optimal transport problem has not been established.

\subsection{Related work}
\paragraph{Convergence of annealed Sinkhorn.} Early convergence results for annealed Sinkhorn were established for discrete optimal assignment problems with uniform marginals on equally sized supports and logarithmic growth of the inverse temperature \citep{sharify_solution_2013}. Potential-preserving state transfer also appears in stabilised scaling algorithms \citep{schmitzer_stabilized_2019}, where stability and iteration bounds are established for temperature transitions initialised from an exact solution at the preceding temperature. More recently, scaling-preserving annealed Sinkhorn was formulated as online mirror descent and analysed under concave inverse-temperature schedules \citep{chizat_annealed_2024}. That analysis, building on \citep{altschuler_near-linear_2017}, identifies an entropic error, controlled by the current temperature, and a relaxation error, controlled by the marginal discrepancy. It provides error bounds for both the actual iterates and a regularisation-path proxy, and motivates a debiased scheme that improves performance in practice. Our work was initially motivated by simplifying the bias correction while retaining its practical benefits, which has ultimately led us back to potential preservation suggested in \citep{schmitzer_stabilized_2019}.

\paragraph{Convergence rates for Sinkhorn.} At fixed temperature Sinkhorn can be viewed as Bregman gradient descent for relative entropy \citep{leger_gradient_2021} and a relative-smoothness analysis gives an $O(t^{-1})$ marginal KL bound \citep{aubin-frankowski_mirror_2022}. Stronger fixed-temperature estimates are also available, including an $O(t^{-2})$ marginal KL convergence rate \cite{ghosal_convergence_2026}, alongside complexity bounds for approximate OT \citep{dvurechensky_computational_2018,altschuler_near-linear_2017}. Geometric convergence is also known for strictly positive kernel matrices, although the classical contraction estimates deteriorate as the temperature tends to zero \citep{schmitzer_stabilized_2019}. Such estimates concern iterations at a fixed temperature and do not directly yield bounds for annealed Sinkhorn, as they depend on the regularisation parameter.

\paragraph{Other Sinkhorn-inspired approaches.} Mirror Sinkhorn combines a mirror-descent step with a projection onto one marginal constraint and provides error bounds for averaged iterates \citep{ballu_mirror_2023}. Annealing appears explicitly in the mirror-descent approach of \citet{kemertas_efficient_2025}, which warm-starts successive entropic optimal transport subproblems and solves each to a prescribed accuracy. The Inexact Proximal point method for exact Optimal Transport (IPOT) uses proximal updates, approximated by warm-started Sinkhorn iterations; its convergence theorem requires control of the errors in these proximal updates, which is not established when only one inner Sinkhorn iteration is performed \citep{xie_fast_2020}. Warm-started IPOT with one inner row–column projection pair has been identified with Bregman Douglas–Rachford splitting and admits an equivalent annealing formulation \citep{chin_annealed_2026}. This iteration differs from potential-preserving annealed Sinkhorn, and its global convergence and a priori convergence rates remain open. We discuss the relation between the two iterations in  Remark~\ref{rmk:ipot}.

\subsection{Our contributions} 
In this manuscript we make the following contributions. 
\begin{enumerate}
    \item We make explicit the distinction between scaling-preserving and potential-preserving annealed Sinkhorn and explain heuristically why preserving the scalings gives rise to the relaxation error identified in \citep{chizat_annealed_2024}. We show that the potential-preserving formulation addresses the same issue as the debiased scheme proposed in \citep{chizat_annealed_2024}, but does so in a more direct way.  We also relate IPOT \cite{xie_fast_2020} to potential-preserving annealing, identifying its warm-started potential transfer as potential preservation with an inertial correction.
    \item We prove that for potential-preserving annealed Sinkhorn under a positive nondecreasing schedule satisfying $\beta_t\to\infty$ and
    \(
    \sum_t \frac{1}{\beta_t}=\infty,
    \)
    there exists a subsequence converging to the unique maximum-entropy optimal transport plan and the full sequence of dual values converges to the value of the optimal transport solution.
    Conversely, when
    \(
    \sum_t \frac{1}{\beta_t}<\infty,
    \)
    we construct an explicit counterexample in which the iterates fail to admit any OT-convergent subsequence, showing that the condition is sharp.
    \item We provide an example satisfying the preceding assumptions for which the full sequence fails to converge. We then prove that the additional condition
    \(
    \sum_{t=1}^{\infty}
    (\frac{\beta_t-\beta_{t-1}}{\beta_{t-1}})^2
    <\infty
    \)
    is sufficient to guarantee convergence of the entire sequence. This condition includes linear inverse schedules.
    \item For polynomial schedules $\beta_t=(t+1)^\kappa$ with \(0 < \kappa \leq 1\), we prove the marginal convergence bound 
    \(
    {\KL(\pi_t\one\mid p) = O\left(t^{-1}\right)}.
    \) For the transport cost after rounding to a feasible plan $R(\pi_t)$, which was introduced in \cite{altschuler_near-linear_2017} (and will be recalled in \eqref{eq:rounding_the_plan}), we also establish ${\langle C,R(\pi_t)\rangle-\OT(p,q)
    =O\left(t^{-1/2}+t^{-\kappa}\right)}$.
    For $1/2\leq\kappa\leq1$, this gives an $O(t^{-1/2})$ bound, improving on the $O(t^{-1/3})$ cost bound for the scaling-preserving iterates and matching the $O(t^{-1/2})$ bound for the regularisation-path proxy of the scaling-preserving scheme~\citep{chizat_annealed_2024}.
\end{enumerate}
Figures~\ref{fig:gap} and~\ref{fig:linear_schedule} provide simple
illustrations of transport cost decay and the role of the annealing
schedule in marginal convergence.

\subsection{Organisation of the paper}
The paper is organised as follows. Section \ref{sec:background} introduces the notation and necessary background on optimal transport and Sinkhorn's algorithm. Section~\ref{subsec:AS} introduces potential-preserving annealed Sinkhorn (PAS) and scaling-preserving annealed Sinkhorn (SAS), explains their differences, and discusses their relation to the debiased scheme \citep{chizat_annealed_2024} and IPOT \citep{xie_fast_2020}. Section \ref{subsec:Conv_theory} presents the main convergence results. For positive, non-increasing regularisation schedules tending to zero, Theorem~\ref{thm:positive} establishes subsequential convergence of the transport plans and convergence of the dual values under the non-summability condition. Proposition~\ref{prop:negative} shows that this condition is sharp, while Example~\ref{ex:non_conv_example} demonstrates that these assumptions do not guarantee convergence of the entire sequence of transport plans. Theorem~\ref{thm:full_convergence} establishes such convergence under an additional square-summability condition on the relative inverse-temperature increments. Theorem \ref{thm:convergence_rate} in Section \ref{subsec:conv_rates} provides a convergence rate bound for PAS, which improves on the rate obtained for SAS in \citep{chizat_annealed_2024}. Section~\ref{sec:Conclusion} summarises the results.
Section~\ref{sec:proofs} contains the proofs.
Appendix~\ref{app:libraries} examines state transfer in computational
optimal transport libraries and the implications of our results
for geometric piecewise-constant annealing schedules.

\section{Background}
\label{sec:background}

\paragraph{Notation.}
We write \(\odot\) and \(\oslash\) for entrywise multiplication and division, \((\cdot)^\top\) for transpose, and \(\operatorname{diag}(u)\) for the diagonal matrix with diagonal \(u\). The vector \(\one\) has all entries equal to one, with its dimension determined by context. Logarithms, exponentials, and powers of vectors or matrices act entrywise. 
We use the Frobenius inner product for matrices.
For nonnegative vectors $x,y$ of the same dimension, we define
$\KL(x\mid y)=\sum_i\bigl(x_i\log(x_i/y_i)-x_i+y_i\bigr)$,
with the conventions $0\log(0/y)=0$ for $y\geq0$ and
$x\log(x/0)=+\infty$ for $x>0$; for matrices, the sum is over
both indices. For vectors, we write $\|x\|_1=\sum_i|x_i|$.
We further write \(\mathbb{R}_+^d=[0,\infty)^d\), \(\mathbb{R}_{++}^d=(0,\infty)^d\), and \(\Delta_d^*=\{r\in\mathbb{R}_{++}^d:r^\top\mathbf1=1\}\). The indicator \(\iota_{\{q\}}(z)\) is zero if \(z=q\) and \(+\infty\) otherwise. The cost oscillation is \(\|C\|_{\mathrm{osc}}=\max_{ij}C_{ij}-\min_{ij}C_{ij}\).

\paragraph{Optimal transport.}
Let \(p\in\Delta_m^*\), \(q\in\Delta_n^*\), and \(C\in\mathbb{R}^{m\times n}_+\) be a cost matrix. A coupling \(\pi\in\mathbb{R}^{m\times n}_+\) specifies the mass transported between each pair of points. The optimal transport problem is
\begin{equation}
    \operatorname{OT}(p,q)
    =\min_{\pi\in\Gamma(p,q)}\langle C,\pi\rangle,
    \qquad
    \Gamma(p,q)
    =\{\pi\geq0:\pi\mathbf 1=p,\ \pi^\top\mathbf 1=q\},
    \label{eq:ot}
\end{equation}
Its entropic regularisation is the strictly convex problem
\begin{equation}
    \operatorname{EOT}_{\beta}(p,q)
    =\min_{\pi\in\Gamma(p,q)}
    \left\{\langle C,\pi\rangle
    +\frac{1}{\beta}\KL(\pi\mid\pi^{\rm ref})\right\},
    \label{eq:eot}
\end{equation}
where \(\beta>0\) is the inverse temperature and \(\pi^{\rm ref}=(mn)^{-1}\mathbf1\mathbf1^\top\) is the uniform reference measure. 

\paragraph{Sinkhorn's algorithm.}
The unique minimiser of \eqref{eq:eot} factorises as \(\pi_\beta = \operatorname{diag}(u)K_\beta\operatorname{diag}(v)\), with \(K_\beta=e^{-\beta C}\) entrywise \citep{peyre_computational_2019,nutz_introduction_2021}. We call \(u\in \R^m_{++}, v\in \R^n_{++}\) the Sinkhorn \emph{scalings}.
Sinkhorn's algorithm iteratively rescales rows and columns to solve the \(\EOT\) problem \citep{cuturi_sinkhorn_2013}:
\begin{equation}
    u\leftarrow p\oslash(K_\beta v),
    \qquad
    v\leftarrow q\oslash(K_\beta^\top u),
    \label{eq:sinkhorn}
\end{equation}
Equivalently, using the \emph{dual potentials}
\(f=\beta^{-1}\log u\) and \(g=\beta^{-1}\log v\), the Sinkhorn updates read
\begin{equation}
    f \leftarrow f+\frac{1}{\beta}\log(p\oslash \pi\one),
    \qquad
    g \leftarrow g+\frac{1}{\beta}\log(q\oslash \pi^\top\one).
    \label{eq:sinkhorn-potentials}
\end{equation}
Stronger regularisation increases bias, while weaker regularisation slows Sinkhorn convergence \citep{ghosal_convergence_2026}.

\section{Annealed Sinkhorn and its convergence}
\label{sec:PAS}
\subsection{Annealing schemes}
\label{subsec:AS}
Annealed Sinkhorn addresses the trade-off in standard Sinkhorn by gradually decreasing \(\varepsilon=\beta^{-1}\) and warm-starting each regularisation level from the current state \citep{schmitzer_stabilized_2019, chizat_annealed_2024}. There are two natural ways to warm-start the algorithm: preserving either the scalings or the dual potentials. To make this distinction precise, fix inverse temperatures \(0<\beta\leq\beta'\), and write
\(
\pi=\operatorname{diag}(u)K_\beta\operatorname{diag}(v)
\) with \(f=\beta^{-1}\log u\) as in Section~\ref{sec:background}. Both schemes first perform the row projection at \(\beta\),
\[
\widehat u=p\oslash(K_\beta v),
\qquad
\widehat f=\beta^{-1}\log\widehat u,
\]
and differ only in the transfer to \(\beta'\).\\
\emph{Scaling-preserving annealed Sinkhorn} (SAS) retains the scaling,
\[
u_{\rm SP}^+:=\widehat u
\; \implies \;
f_{\rm SP}^+=\frac{\beta}{\beta'}\widehat f.
\]
\emph{Potential-preserving annealed Sinkhorn} (PAS) retains the
potential instead,
\[
f_{\rm PP}^+:=\widehat f
\; \implies \;
u_{\rm PP}^+=\widehat u^{\,\beta'/\beta}.
\]

Thus, one complete annealed Sinkhorn step from \((\beta,u,v)\) to \((\beta',u^+,v^+)\) is
\[
\widehat u=p\oslash(K_\beta v),\qquad
u^+=
\begin{cases}
\widehat u, & \text{\rm SAS},\\
\widehat u^{\,\beta'/\beta}, & \text{\rm PAS},
\end{cases}
\qquad
v^+=q\oslash(K_{\beta'}^\top u^+).
\]
The two schemes therefore coincide at fixed temperature, but differ in their transfer rule when \(\beta'>\beta\). In particular, if \(\pi\) solves \({\rm EOT}_\beta(p,q)\), then \(\widehat u=u\) and \(\widehat f=f\), hence
\[
f_{\rm SP}^+=\frac{\beta}{\beta'}f,
\qquad
f_{\rm PP}^+=f.
\]
Thus SAS contracts the dual potential solely because the temperature changes, whereas PAS preserves it. This distinction is also apparent from
\[
\pi_{ij}=u_i e^{-\beta C_{ij}}v_j
=\exp\!\bigl(\beta(f_i+g_j-C_{ij})\bigr).
\]

\begin{figure}
    \centering
    \includegraphics[width=\textwidth]{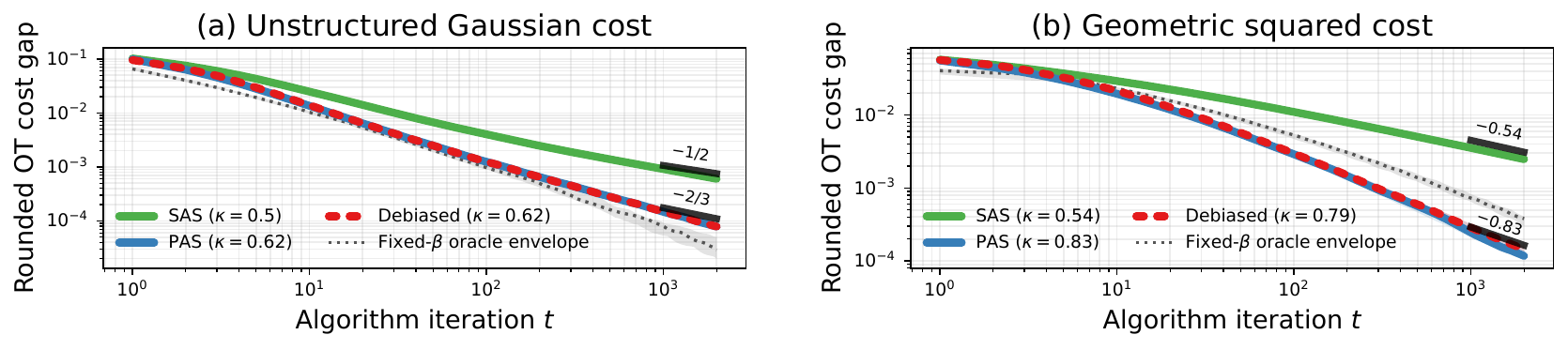}
    \caption{Illustration of the rounded transport cost gap
    $G_t=\langle C,R(\pi_t)\rangle-\OT(p,q)$, where the same rounding
    to the prescribed marginals is applied to every method
    (see~\eqref{eq:rounding_the_plan}).
    (a) A $100\times100$ cost matrix with independent standard Gaussian
    entries and marginals obtained by normalising independent
    $\operatorname{Unif}(0,1)$ weights.
    (b) Squared Euclidean cost between two planar point clouds of
    300 points each, with uniform weights.
    Both cost matrices are shifted and rescaled to $\min C=0$ and
    $\max C=1$.
    Scaling-preserving (SAS), potential-preserving (PAS), and debiased
    annealing use $\beta_t=10(1+t)^\kappa$; the exponents shown
    in the legends are rounded and were selected separately for each
    method and problem class to minimise the median gap at 200
    iterations over $\kappa\in\{0,1/24,\ldots,1\}\cup\{1.1\}$.
    All runs start from $u=v=\one$, and an iteration consists of one
    row and one column update.
    Curves show medians and bands show interquartile ranges over
    20 random instances. The grey dotted curve gives the pointwise
    best result among 25 fixed inverse temperatures spaced
    geometrically from 10 to 4800; black segments indicate reference
    slopes.}
    \label{fig:gap}
\end{figure}

Intuitively, the quantities \(f_i,g_j\) are ``prices'' in the same units as \(C_{ij}\), while \(u,v\) are their exponentiated representations at the current temperature. It is therefore natural for the scalings to vary with the temperature while the potentials remain fixed.

Appendix~\ref{app:libraries} examines the state-transfer rules
used in computational OT libraries. In the following two
remarks, we use this viewpoint to compare PAS with debiased
SAS and IPOT.

\begin{remark}[Relation to debiased SAS]
\label{rmk:debiased}
The debiased SAS scheme in \citep{chizat_annealed_2024}
corrects the scaling-preserving transfer as follows.
Writing $0<\beta^-<\beta<\beta'$ for three consecutive
inverse temperatures, its updates are
\[
u^+=\widehat u\,u^{\,1-\beta^-/\beta}
\qquad
\text{and}
\qquad
f^+
=
\frac{\beta}{\beta'}\widehat f
+
\frac{\beta-\beta^-}{\beta'}f.
\]
At an exact solution for the current inverse temperature \(\beta\), this reduces to
\(
f^+
=
\frac{2\beta-\beta^-}{\beta'}f.
\)
For a linear inverse-temperature schedule, \(\beta'-\beta=\beta-\beta^-\), so at such a solution the debiased transfer is exactly potential-preserving. More generally, it is a lagged correction of the scaling-preserving transfer, whereas PAS directly preserves the new dual potential.
\end{remark}
\begin{remark}[Relation to IPOT]\label{rmk:ipot}
The \emph{Inexact Proximal point method for exact Optimal Transport
(IPOT)} \citep{xie_fast_2020} approximates, at each iteration,
the proximal point
\[
    \widehat P_t
    :=\operatorname*{argmin}_{\pi\in\Gamma(p,q)}
    \bigl\{\langle C,\pi\rangle+\eta\KL(\pi\mid P_{t-1})\bigr\},
\]
where $\eta>0$ is fixed. This is the entropic problem
\eqref{eq:eot} with inverse temperature $\eta^{-1}$ and
reference $P_{t-1}$ in place of $\pi^{\rm ref}$. Its solution
$\widehat P_t$ is a diagonal scaling of
$Q_t=P_{t-1}\odot e^{-C/\eta}$.
IPOT performs finitely many updates \eqref{eq:sinkhorn},
with $K_\beta$ replaced by $Q_t$, and sets
$P_t=\operatorname{diag}(u)Q_t\operatorname{diag}(v)$.
The annealing interpretation follows from
\[
    P_{t-1}
    =\operatorname{diag}(a)e^{-\beta_{t-1}C}
      \operatorname{diag}(b)
    \quad\Longrightarrow\quad
    Q_t
    =\operatorname{diag}(a)
      e^{-(\beta_{t-1}+\eta^{-1})C}
      \operatorname{diag}(b).
\]
With $P_0=\pi^{\rm ref}$, induction gives $\beta_t=t/\eta$.
The kernels $Q_t$ and $K_{\beta_t}$ differ only by positive
diagonal factors. Absorbing these factors into the inner
scaling vectors transforms the IPOT updates into
\eqref{eq:sinkhorn} at inverse temperature $\beta_t$.
Thus each outer iteration consists of finitely many
Sinkhorn updates for \eqref{eq:eot} at the next temperature.
Under this linear schedule, SAS converges to a
semi-unbalanced OT solution
\citep[Theorem~2.1]{chizat_annealed_2024}.
IPOT, however, retains the inner scalings between outer
iterations, which differs from preserving the accumulated
scalings $a,b$. For one row and one column update per outer
iteration, an equivalent recurrence for the accumulated
scalings is derived in \citep[Proposition~1]{chin_annealed_2026}.
In terms of the dual potentials, this recurrence gives
\[
    \widetilde g
    =\frac{2\beta}{\beta'}g-\frac{\beta^-}{\beta'}g^-
    =g+\frac{\beta^-}{\beta'}(g-g^-),
    \qquad \beta'=2\beta-\beta^-,
\]
where $g^-,g$ are the column potentials at consecutive
inverse temperatures $\beta^-,\beta$.
Thus IPOT adds to $g$ a multiple of the preceding increment
$g-g^-$, an inertial correction to the potential-preserving
transfer. This relates it to PAS, which converges to an
OT solution under linear annealing
(Theorem~\ref{thm:full_convergence}).
The inertial correction and the use of the new temperature
for both marginal projections prevent a direct application
of our convergence results to IPOT.
\end{remark}

\subsection{Convergence results}
\label{subsec:Conv_theory}
We now fix the notation for repeated PAS updates along an
inverse-temperature schedule $(\beta_t)_{t\ge0}$.
The PAS iterates are defined by the following algorithm.
\begin{algorithm}[H]
\caption{Potential-preserving Annealed Sinkhorn (PAS).}
\label{alg:PAS}

\begin{enumerate}
    \item \textbf{Input:}
    probability vectors
    $p \in \Delta_m^*$,
    $q \in \Delta_n^*$,
    cost $C \in \mathbb{R}_+^{m \times n}$,
    positive and nondecreasing annealing schedule $(\beta_t)_{t \geq 0}$

    \item \textbf{Initialize:}
    let $u_0 = \mathbf{1} \in \mathbb{R}^m$,
    $v_0 = \mathbf{1} \in \mathbb{R}^n$,
    $K_0 = e^{-\beta_0 C} \in \mathbb{R}_{++}^{m \times n}$,
    and $\pi_0 = \operatorname{diag}(u_0)K_0\operatorname{diag}(v_0)=K_0$

    \item \textbf{For} $t = 1,2,\ldots$ let
    \begin{align*}
        u_t &= (p \oslash (K_{t-1}v_{t-1}))^\frac{\beta_t}{\beta_{t-1}},
        && \text{\# project on row marginal constraint and rescale}\\
        K_t &= e^{-\beta_t C},
        && \text{\# update inverse temperature}\\
        v_t &= q \oslash (K_t^\top u_t),
        && \text{\# project on column marginal constraint}\\
        \pi_t &= \operatorname{diag}(u_t)
        K_t \operatorname{diag}(v_t)
        && \text{\# define primal iterate}
    \end{align*}
\end{enumerate}

\end{algorithm}
The complete iterates $(\pi_t)_{t\ge1}$ belong to
\begin{equation}\label{eq:Gamma_star_q}
    \Gamma(*,q)
    :=\{\pi\in\mathbb R_+^{m\times n}:\pi^\top\one=q\}.
\end{equation}
Since $q$ has unit mass, these iterates are probability matrices.

Our convergence analysis builds on the online mirror-descent
viewpoint of \citep{chizat_annealed_2024,ballu_mirror_2023},
adapted here to potential-preserving annealing.
To identify the limiting plan, we first introduce the
maximum-entropy optimal transport solution. Among the solutions of the optimal transport problem \eqref{eq:ot} given by \(\mathcal S_{\rm OT}\),
let $\pi^\star$ denote the unique plan that minimises the negative entropy
\[
\label{C1}\tag{MaxEnt}
\pi^\star
\coloneq
\operatorname*{arg\,min}_{\gamma\in\mathcal S_{\rm OT}}
\mathcal H(\gamma),
\qquad
\mathcal H(\pi):=\sum_{ij}\pi_{ij}\log\pi_{ij},
\qquad
\mathcal S_{\rm OT}
:=\operatorname*{argmin}_{\pi\in\Gamma(p,q)}\langle C,\pi\rangle
\]
with $0\log0=0$.
We call $\pi^\star$ the maximum-entropy optimal transport plan. The minimiser in \eqref{C1} is unique since $\mathcal H$ is strictly convex and the set of optimal plans \(\mathcal S_{\rm OT}\) is convex. We first prove the following.
\begin{theorem}[Subsequential convergence and convergence of the dual value]
\label{thm:positive}
Let $(\beta_t)_{t\ge0}$ be positive and nondecreasing, with
\begin{equation}\label{AS1}\tag{A1}
    \beta_t\to\infty,
    \qquad
    \sum_{t=0}^{\infty}\beta_t^{-1}=\infty.
\end{equation}
Let $(u_t, v_t)_{t\geq 0}$ and $(\pi_t)_{t \geq 0}$ be the corresponding PAS scalings and plans generated by Algorithm \ref{alg:PAS} and let $\pi^\star$ denote the maximum-entropy optimal transport plan defined in \eqref{C1}. Then:
\begin{enumerate}
\item there exists a subsequence $(t_j)$ such that $\pi_{t_j}\longrightarrow\pi^\star$;
\item the dual values
\[
L_t:=p^\top f_t+q^\top g_t,
\qquad
f_t:=\beta_t^{-1}\log u_t,
\quad
g_t:=\beta_t^{-1}\log v_t,
\]
satisfy $L_t\uparrow\OT(p,q)$ as $t\to\infty$, with $(L_t)_{t\ge1}$ nondecreasing.
\end{enumerate}
\end{theorem}
\begin{proof}[Proof sketch]
The online mirror-descent representation of PAS yields
a one-step KL identity whose summation gives, for $1\le a<T$,
\begin{equation}\label{eq:PAS_telescope}
\begin{aligned}
\sum_{t=a+1}^{T}
\beta_{t-1}^{-1}\KL(\pi_t\one\mid p)
\le{}
\beta_a^{-1}\KL(\pi^\star\mid\pi_a)
+(\beta_a^{-1}-\beta_T^{-1})\log(mn).
\end{aligned}
\end{equation}
For fixed $a$, the right-hand side is bounded independently
of $T$. Since $\sum_t\beta_t^{-1}=\infty$, the marginal
discrepancy $\KL(\pi_t\one\mid p)$ has liminf zero. Along a subsequence on which this discrepancy vanishes, Pinsker's inequality gives $\pi_t\one\to p$. Since $\pi_t^\top\one=q$, compactness provides a further convergent subsequence, whose limit is a solution to the OT problem by the standard stability result recalled in Lemma~\ref{lem:B1_Chizat}. An entropy comparison between optimal plans, using the Gibbs representation of $\pi_t$, identifies this limit as $\pi^\star$. 
The same identity also shows that $(L_t)$ is nondecreasing,
while the Gibbs representation gives $L_t\to\OT(p,q)$ along the convergent subsequence, hence
along the entire sequence. The full proof is given in Section~\ref{app:convergence}.
\end{proof}

The following proposition shows that the non-summability assumption in \eqref{AS1} cannot be omitted: for every summable temperature schedule, there exists an OT problem for which no subsequence converges to an optimal transport plan.
\begin{proposition}
\label{prop:negative}
Let $(\beta_t)_{t\ge0}$ be a positive, nondecreasing
inverse-temperature schedule satisfying
\[
\sum_{t=0}^{\infty}\beta_t^{-1}<\infty.
\]
Then there exist a cost matrix $C\in\mathbb R_+^{2\times2}$
and probability vectors $p,q\in\Delta_2^*$ such that the
corresponding PAS iterates of Algorithm \ref{alg:PAS} satisfy
\[
\pi_t\longrightarrow\pi_\infty,
\qquad
\pi_\infty\one\neq p.
\]
\end{proposition}
The proof is given in Section~\ref{app:proof_prop}.

Condition \eqref{AS1} alone does not imply convergence of the entire sequence. 
\begin{example}[Failure of full convergence]
\label{ex:non_conv_example}
Let $p=(a,1-a)^\top$, where
$a\in(0,1)\setminus\{1/2\}$, and take $q=(1)$ and
$C=(0,1)^\top$. Set $\beta_0=1$ and
\[
\beta_t=2^k
\qquad\text{for }2^k\le t<2^{k+1},\quad k\ge0.
\]
This inverse-temperature schedule satisfies \eqref{AS1},
but the first marginal $\pi_t\one$ has two distinct
accumulation points, while $L_t\uparrow\OT(p,q)=1-a$.
\end{example}

The calculation is given in
Section~\ref{sec:example_calculation}.

Example~\ref{ex:non_conv_example} shows that abrupt changes in
temperature can prevent full convergence even when \eqref{AS1}
holds. Controlling these changes allows us to strengthen
Theorem~\ref{thm:positive}: if the relative increments of the
inverse temperature are square-summable, the entire sequence
converges to the maximum-entropy optimal transport plan.

\begin{theorem}[Full convergence]
\label{thm:full_convergence}
Let $(\pi_t)$ be the PAS iterates of Algorithm \ref{alg:PAS}
for a positive, nondecreasing inverse-temperature schedule
$(\beta_t)_{t\ge0}$ satisfying \eqref{AS1} and let \(\pi^\star\) be the unique maximum-entropy transport plan given by \eqref{C1}. If, in addition,
\begin{equation}\label{AS2}\tag{A2}
    \sum_{t=1}^{\infty}
    \left(\frac{\beta_t-\beta_{t-1}}{\beta_{t-1}}\right)^2
    <\infty,
\end{equation}
then
\[
\lim_{t\to\infty}\KL(\pi^\star\mid\pi_t)=0.
\]
\end{theorem}
\begin{proof}[Proof sketch]
The main step is to establish the estimate
\begin{equation}\label{eq:marginal_step}
\KL(\pi_t\one\mid p)
\le
\KL(\pi_{t-1}\one\mid p)
-\frac12\KL(\pi_{t-1}\mid\pi_t)
+C_0\left(\frac{\beta_t-\beta_{t-1}}{\beta_t}\right)^2,
\end{equation}
for $t\ge2$, where $C_0>0$ is independent of $t$. We obtain this by
expressing the PAS update as a row projection followed by
a change of temperature and a column projection. A uniform
logarithmic-moment bound controls the temperature-change
term, which is then absorbed into the KL decrease using
Young's inequality. Under \eqref{AS2}, the remaining error
is summable, so $\KL(\pi_t\one\mid p)$ converges. By 
Theorem~\ref{thm:positive}, this limit is zero. Every accumulation point therefore solves \eqref{eq:ot} and the entropy-selection argument identifies
it with $\pi^\star$. Compactness then gives convergence
of the entire sequence, and consequently
$\KL(\pi^\star\mid\pi_t)\to0$.
The full proof is given in
Section~\ref{app:thm_full_conv_proof}.
\end{proof}

\begin{remark}[Schedule assumptions]
\label{rmk:guarantee}
For a positive, nondecreasing schedule with $\beta_t\to\infty$,
the natural bounded increment condition
\[
\sup_{t\ge1}(\beta_t-\beta_{t-1})<\infty
\]
implies both \eqref{AS1} and \eqref{AS2}. In particular, polynomial schedules $\beta_t=(1+t)^\kappa$
with $0<\kappa\le1$ satisfy this bound, so Theorem~\ref{thm:full_convergence} covers all such schedules, including linear annealing. For the linear schedule $\beta_t=1+t$, this contrasts with SAS, whose iterates converge to a solution of a semi-unbalanced OT problem \citep{chizat_annealed_2024}. Conditions \eqref{AS1} and \eqref{AS2} also admit schedules with unbounded increments, such as $\beta_t=(t+2)\log(t+2)$. Further iterated logarithmic factors may be included, with the shift chosen sufficiently large to keep all factors positive. In contrast, the doubling schedule in Example~\ref{ex:non_conv_example} satisfies \eqref{AS1}
but violates \eqref{AS2}.
\end{remark}

\begin{figure}
    \centering
    \includegraphics[width=\textwidth]{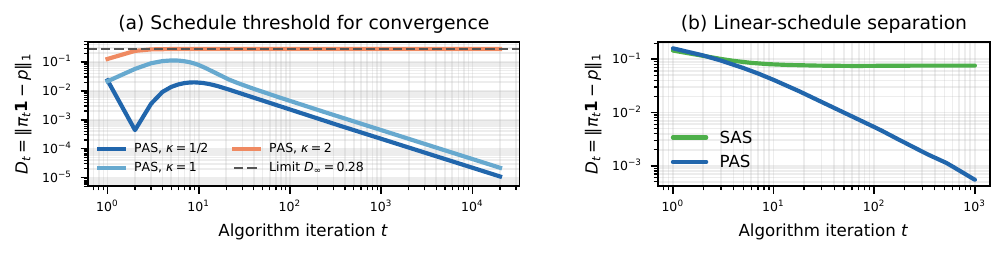}
    \caption{Illustration of the schedule conditions for convergence,
    measured by the source marginal error
    $D_t=\|\pi_t\one-p\|_1$; the target marginal is exactly $q$ at
    every iteration. All runs start from $u=v=\one$.
    (a) PAS for $p=(0.64,0.36)$, $q=(0.5,0.5)$ and the $2\times2$
    cost matrix with zero diagonal and unit off-diagonal entries,
    under $\beta_t=(1+t)^\kappa$.
    For $\kappa\in\{1/2,1\}$, Theorem~\ref{thm:full_convergence}
    gives $D_t\to0$. For $\kappa=2$, the construction in
    Proposition~\ref{prop:negative} gives
    $\pi_t\to\operatorname{diag}(1/2,1/2)$ and $D_t\to0.28$
    (grey dashed line).
    (b) PAS and SAS for 96 source and 96 target points sampled from
    two-component Gaussian mixtures in $\mathbb R^2$, with normalised
    log-normal weights (log standard deviation $0.7$) and squared
    Euclidean cost divided by its median, under $\beta_t=1+0.2t$.
    The decreasing PAS error and the nonzero SAS plateau illustrate
    the distinction between balanced and semi-unbalanced limits
    discussed in Remark~\ref{rmk:guarantee}.}
    \label{fig:linear_schedule}
\end{figure}

\section{Convergence rates}
\label{subsec:conv_rates}

For polynomial inverse-temperature schedules, the estimates
underlying Theorems~\ref{thm:positive} and
\ref{thm:full_convergence} also yield quantitative bounds
on the marginal error and the transport cost.
Since the PAS iterates need not satisfy both marginal
constraints, we evaluate the cost after rounding:
let $R(\pi_t)\in\Gamma(p,q)$ denote the feasible plan
obtained by the rounding procedure of
\citep{altschuler_near-linear_2017}, which, for $\pi \in \Gamma(*,q)$ (defined in \eqref{eq:Gamma_star_q}), is given by
\begin{equation}\label{eq:rounding_the_plan}
\widehat\pi
=
\operatorname{diag}\!\left(
\min\{\mathbf1,p\oslash(\pi\mathbf1)\}
\right)\pi,
\qquad
R(\pi)
:=
\widehat\pi+
\frac{(p-\widehat\pi\mathbf1)
(q-\widehat\pi^\top\mathbf1)^\top}
{1-\mathbf1^\top\widehat\pi\mathbf1}.
\end{equation}
If $\pi\one=p$, we simply set $R(\pi)=\pi$.

\begin{theorem}[Convergence rates]
\label{thm:convergence_rate}
Let $(\pi_t)$ be the PAS iterates of Algorithm \ref{alg:PAS}
for the inverse-temperature schedule
$\beta_t=(t+1)^\kappa$, where $0<\kappa\le1$.
There exists a constant $A>0$, depending on $p,q,C$ and
$\kappa$ but independent of $t$, such that
\[
\KL(\pi_t\one\mid p)\le A t^{-1},
\qquad t\ge1.
\]
Consequently, Pinsker's inequality and the rounding estimate recalled in \eqref{eq:rounding_the_plan} give
\[
\|\pi_t\one-p\|_1\le\sqrt{2A}\,t^{-1/2},
\qquad
0\le\langle C,R(\pi_t)\rangle-\OT(p,q)
\le B\bigl(t^{-\kappa}+t^{-1/2}\bigr),
\]
for every $t\ge1$, where $B>0$ depends only on $p,q,C$
and $\kappa$.
\end{theorem}

\begin{proof}[Proof sketch]
Theorem~\ref{thm:full_convergence} gives a uniform bound on
$\KL(\pi^\star\mid\pi_t)$. Applying the telescoping estimate \eqref{eq:PAS_telescope} from the proof of Theorem~\ref{thm:positive} over $\lfloor T/2\rfloor<t\le T$ shows that some iterate in this interval has marginal KL discrepancy bounded by a constant multiple of $T^{-1}$. The one-step estimate \eqref{eq:marginal_step} from the proof of Theorem~\ref{thm:full_convergence} controls its subsequent
increase: since $(\beta_t-\beta_{t-1})/\beta_t\le\kappa t^{-1}$,
the accumulated error over the remaining iterations is
also bounded by a constant multiple of $T^{-1}$. This gives the claimed bound at iteration $T$. Pinsker's inequality gives the marginal $\ell^1$ estimate and combining it with the rounding and entropic-bias estimates of \citep[Lemma~2.3]{chizat_annealed_2024} gives the cost bound.
\end{proof}

\begin{remark}[Comparison with fixed-temperature Sinkhorn]
Theorem~\ref{thm:convergence_rate} extends the $\mathcal{O}(t^{-1})$ marginal KL bound known for standard Sinkhorn ($\kappa=0$) \citep[Proposition~8]{aubin-frankowski_mirror_2022}
to potential-preserving annealing with $0<\kappa\le1$.
Faster rates are known at fixed temperature, including $\mathcal{O}(t^{-2})$ marginal KL bounds under suitable assumptions \citep[Corollary 4.6]{ghosal_convergence_2026}.
\end{remark}
\begin{remark}[Comparison with SAS]
For $1/2\le\kappa\le1$, Theorem~\ref{thm:convergence_rate}
gives an $O(t^{-1/2})$ bound on the rounded cost gap of
the PAS iterates. This improves on the $O(t^{-1/3})$
bound obtained for the SAS iterates in
\citep[Section~2.3]{chizat_annealed_2024} and matches the
$O(t^{-1/2})$ bound established there for the
regularisation-path proxy. The PAS bound need not be sharp. Its decay exponent is the same throughout $\kappa\in[1/2,1]$, although the constants may depend on $\kappa$. Sharper estimates could determine whether the decay exponent
depends on $\kappa$.
\end{remark}

\section{Conclusion}
\label{sec:Conclusion}

We established a convergence theory for potential-preserving annealed Sinkhorn, including a sharp criterion for subsequential convergence, sufficient conditions for convergence of the entire sequence, and quantitative bounds on the marginal and transport cost errors.
Through an online mirror-descent formulation, we proved that, for positive nondecreasing inverse temperatures tending to infinity, the non-summability condition $\sum_t\beta_t^{-1}=\infty$ implies subsequential convergence to the maximum-entropy optimal transport plan and convergence of the dual objective to $\OT(p,q)$.
For every summable temperature schedule, we constructed an explicit $2\times2$ counterexample in which the iterates converge to a plan that violates the source marginal constraint. Thus the non-summability condition is sharp: without it, even subsequential convergence to an
optimal transport plan may fail.
We also showed that convergence of the dual values does not imply convergence of the transport plans, by constructing an example with repeated temperature jumps. Under the additional assumption \(\sum_{t=1}^{\infty}
((\beta_t-\beta_{t-1})/\beta_{t-1})^2<\infty\), we established convergence of the entire sequence to the maximum-entropy optimal transport plan. In particular, this includes linear inverse-temperature schedules, under which SAS converges instead to a solution of a semi-unbalanced OT problem.
For polynomial schedules $\beta_t=(t+1)^\kappa$ with $0<\kappa\leq1$, we obtained marginal KL and $\ell^1$ bounds of $O(t^{-1})$ and $O(t^{-1/2})$, respectively, together with a rounded OT cost gap bound of $O(t^{-\kappa}+t^{-1/2})$. 
Figures~\ref{fig:gap} and~\ref{fig:linear_schedule} illustrate
the decay of the rounded transport cost gap and the dependence
of marginal convergence on the annealing schedule.
Appendix~\ref{app:libraries} discusses state-transfer conventions in computational OT libraries and applies the subsequential convergence criterion to geometric piecewise-constant schedules.

\paragraph{Acknowledgements.} Both authors are supported by Research England through the Expanding Excellence in England (E3) award to MARS: Mathematics for AI in Real-world Systems, in the School of Mathematical Sciences at Lancaster University. We thank Jean Feydy for discussions that prompted us to examine the connection between
IPOT and annealing, and for drawing our attention to the prior work on annealed Sinkhorn \citep{chizat_annealed_2024}. These discussions took place during Dagstuhl Seminar 25492, \emph{Generalized Voronoi Diagrams and Applications}, held at Schloss Dagstuhl in December 2025. 

\paragraph{Research provenance and use of AI.}
The work began in summer 2025, when we considered a
fixed-point approach to solving optimal transport problems, provisionally
called ``Fixhorn'', in which each Sinkhorn iterate served
as the reference measure for the next iteration. We subsequently recognised a connection with IPOT
\citep{xie_fast_2020}. Discussions with Jean Feydy at
Dagstuhl Seminar 25492 then led us to examine its
relation to annealing and to study
\citep{chizat_annealed_2024}.
The convergence of SAS to a semi-unbalanced OT solution under linear inverse-temperature schedules led us to investigate whether modifying the source measure could recover convergence to a solution of the prescribed balanced OT problem. Comparing the optimality conditions suggested a correction depending on the unknown optimal potentials of the associated balanced entropic OT problem, which we sought to approximate using the current iterates. This idea, akin to the marginal reweighting
in \citep{liu_solving_2025} but reached independently, motivated our initial modification: raising the row scaling
$p\oslash(K_{t-1}v_{t-1})$ to the power
$\theta_t=(\alpha_t+\beta_t)/\beta_{t-1}$, where
$\alpha_t=\beta_t-\beta_{t-1}$.
We sketched an online
mirror-descent formulation for what we then regarded as an
accelerated or compensated annealed Sinkhorn method.
These calculations led us to propose simplifying the exponent from $\theta_t$ to the ratio $\beta_t/\beta_{t-1}$ used in Algorithm~\ref{alg:PAS}, which simplified the convergence argument and led to an initial proof of Theorem~\ref{thm:positive}. We soon recognised that the modified update preserved the dual potentials and connected it with \citep{schmitzer_stabilized_2019}. Up to this point, the work progressed largely without AI assistance. We subsequently used GPT-5.6 and GPT-6 models to help simplify the proof of Theorem~\ref{thm:positive}, develop Theorems~\ref{thm:full_convergence} and~\ref{thm:convergence_rate}, and 
construct the counterexamples in Proposition~\ref{prop:negative} and Example~\ref{ex:non_conv_example}. Specific questions about extending our subsequential result led the models to suggest an estimate of the form
\[
    d_{t+1}\le d_t-cS_t+e_t,
    \qquad d_t,S_t,e_t\ge0,\quad c>0,\quad
    \sum_t e_t<\infty.
\]
Here $d_t$ is the marginal KL discrepancy and $S_t$ measures the change between successive iterates. Such an estimate implies convergence of $d_t$, so the previously established $\liminf_t d_t=0$ gives $d_t\to0$.
This is a deterministic instance of the Robbins--Siegmund convergence theorem \citep{robbins_convergence_1971} and uses the same summable-error principle as quasi-Fej\'er analysis \citep{combettes_quasi_2001}. 
All mathematical arguments were checked by the authors and we take full responsibility for the results.

\section{Proofs} 
\label{sec:proofs}

Throughout the proofs, in addition to the notation introduced in Section~\ref{sec:background}, we use the following conventions
\begin{align}
\label{eq:notation}
    \alpha_t\coloneq \beta_t-\beta_{t-1}, \qquad  \theta_t \coloneq \frac{\beta_t}{\beta_{t-1}}, \qquad \mathcal H(\pi)=\sum_{ij}\pi_{ij}\log\pi_{ij}.
\end{align}
We will use the notation \(\pi\in\Gamma(p,q)\) for an arbitrary transport plan, \(\bar\pi\in\operatorname*{argmin}_{\pi\in\Gamma(p,q)}\langle C,\pi\rangle\) for an optimal transport plan, and $\pi^\star$ for the unique maximum-entropy transport plan, as defined in
\eqref{C1}.
We also utilise key properties of KL divergence which can be found in~\citep{chen_convergence_1993,cover_elements_2006}, as well as the convergence result for entropic optimal transport from~\citep[Lemma B.1]{chizat_annealed_2024} given by the following Lemmas.

\begin{lemma}[Standard KL properties]
\label{lem:KL_props}
Let $X,Y,Z\in\R_+^{m\times n}$, with $\KL$ defined using
the conventions in Section~\ref{sec:background}.
Then the following properties hold.
\begin{enumerate}
\item Non-negativity
\begin{align}
\label{eq:KL_non_neg}
    \KL(X\mid Y) \geq 0.
\end{align}

\item Data processing inequality
\begin{align}
\label{eq:data_processing_ineq}
    \KL(X\mid Y)\ge \KL(X\one\mid Y\one).
\end{align}

\item Three-point identity. If $X,Y\in\R_{++}^{m\times n}$, then
\begin{align}
\label{eq:three_point_identity}
    \left\langle
    \log(X\oslash Y)
    ,X-Z
    \right\rangle
    =
    \KL(X\mid Y)
    +\KL(Z\mid X)
    -\KL(Z\mid Y).
\end{align}

\item Pinsker's inequality.
If \(\sum_{ij}X_{ij} = \sum_{ij}Y_{ij} = 1\) then
\begin{align}
\label{eq:pinskers_ineq}
2\KL(X\mid Y) \geq \|X-Y\|_1^2.
\end{align}
\end{enumerate}
\end{lemma}
Note that since we use the generalised definition of KL divergence, \eqref{eq:KL_non_neg} holds irrespective of the total masses of $X$ and $Y$.
\begin{lemma}[Lemma B.1 in \citep{chizat_annealed_2024}]
\label{lem:B1_Chizat}
For any $\beta_t>0$, $u_t\in(\mathbb{R}_{++})^m$, and $v_t\in(\mathbb{R}_{++})^n$, let ${\pi_t=\operatorname{diag}(u_t)e^{-\beta_t C}\operatorname{diag}(v_t)}$ and consider the marginals $p_t=\pi_t\one\in\Delta_m^*$ and $q_t=\pi_t^\top\one\in\Delta_n^*$. Then $\pi_t$ is the unique minimizer for $\mathrm{EOT}_{\beta_t}(p_t,q_t)$, and $(f_t,g_t)\coloneqq(\beta_t^{-1}\log(u_t/n),\beta_t^{-1}\log(v_t/m))$ is a maximizer of the dual EOT problem.

Considering now sequences indexed by $t$, if we moreover have $\beta_t\to\infty$, $p_t\to p$ and $q_t\to q$ then $\mathrm{EOT}_{\beta_t}(p_t,q_t)\to\mathrm{OT}(p,q)$ and, up to taking a subsequence, $\pi_t$ converges to an optimal transport plan $\bar\pi$ and $(f_t,g_t)$ converge in the norm $\|\cdot\|_{o,\infty}$ to Kantorovich potentials, i.e.\ solutions to the dual of $\mathrm{OT}(p,q)$.
\end{lemma}

\subsection{Proof of Theorem~\ref{thm:positive}}
\label{app:convergence}

We will first establish the online mirror-descent formulation of PAS and derive an identity for the change in
$\beta_t^{-1}\KL(\pi^\star\mid\pi_t)$. We will then prove
subsequential convergence to an optimal transport plan
and convergence of the dual values. Finally, we will identify the limiting plan as $\pi^\star$.

All one-step identities below are used for $t\ge2$, so that
$\pi_{t-1}$ and $\pi_t$ both have second marginal $q$ and
unit mass. Let \(K_t \coloneq K_{\beta_t} = e^{-\beta_t C}\) and define the positive reference vector
\(
    p_t=p\odot u_{t-1}^{(\theta_t-1)/\theta_t}.
\) The vector $p_t$ depends on the preceding PAS iterate and it need not have
unit mass. The prescribed source marginal of the algorithm remains $p$.

\subsubsection{Online mirror-descent formulation}

We show that the PAS update from $\pi_{t-1}$ to $\pi_t$
in Algorithm~\ref{alg:PAS} is an online mirror-descent (OMD) step. 
For the OMD formulation of scaling-preserving annealed Sinkhorn and the underlying OMD framework, see
\citep[Section~2.2 and Appendix~C]{chizat_annealed_2024}.
Consider the following OMD step from \(\pi_{t-1}\):
\begin{equation*}
    \tilde{\pi}_t
    =
    \argmin_{\pi\in\mathbb{R}_+^{m\times n}}
    \left\{
        \langle \alpha_t C,\pi\rangle
        +
        \left\langle
            \theta_t
            \log\left(
                (\pi_{t-1}\one) \oslash p_t
            \right)
            \one^\top,
            \pi
        \right\rangle
        +
        \iota_{\{q\}}(\pi^\top\one)
        +
        \KL(\pi\mid\pi_{t-1})
    \right\}.
\end{equation*}
We show that \(\widetilde\pi_t=\pi_t\), where \(\pi_t\) is given by PAS. To characterise this update, introduce a Lagrange multiplier
\(\lambda\in\mathbb{R}^n\) associated with the constraint
\(\pi^\top\one=q\), which is enforced through the indicator function \(\iota_{\{q\}}(\pi^\top\one)\) in the equation above. The corresponding Lagrangian is
\begin{equation}
    \label{eq:OMG_Lagrangian}
    L(\pi,\lambda)
    \coloneqq\;
        \langle \alpha_t C,\pi\rangle
        +
        \left\langle
            \theta_t
            \log\left(
                (\pi_{t-1}\one) \oslash p_t
            \right)
            \one^\top,
            \pi
        \right\rangle
        +
        \langle \lambda,\pi^\top\one-q\rangle
        +
        \KL(\pi\mid\pi_{t-1}).
\end{equation}
The associated dual problem is therefore
\begin{equation*}
    \max_{\lambda} g(\lambda),
    \qquad
    g(\lambda)
    \coloneqq
    \min_{\pi\in\mathbb{R}_+^{m\times n}} L(\pi,\lambda).
\end{equation*}
Since the objective is strictly convex, it suffices to verify the first-order optimality condition at the feasible point $\pi_t$. Namely, to find $\lambda$ such that
\[
    0=
    \nabla_\pi L(\pi,\lambda)\big|_{\pi=\pi_t}=
    \alpha_t C
    +
    \theta_t
    \log\left(
        (\pi_{t-1}\one) \oslash p_t
    \right)
    \one^\top
    +
    \one\lambda^\top
    +
    \log(\pi_t\oslash\pi_{t-1}),
\]
with \(\pi_t^\top\one=q\).
Write the consecutive iterates of PAS in the standard form as
\(
    \pi_{t-1}
    =
    \operatorname{diag}(u_{t-1})
    K_{t-1}
    \operatorname{diag}(v_{t-1}),
\)
and
\(
    \pi_t
    =
    \operatorname{diag}(u_t)
    K_t
    \operatorname{diag}(v_t).
\)
The first marginal of \(\pi_{t-1}\) is therefore
\(
    \pi_{t-1}\one
    =
    u_{t-1}\odot(K_{t-1}v_{t-1}),
\)
while
\[
    \log(\pi_t\oslash\pi_{t-1})
    =
    \log(K_t\oslash K_{t-1})
    +
    \log(u_t\oslash u_{t-1})\one^\top
    +
    \one
    \log(v_t\oslash v_{t-1})^\top.
\]
Substituting these expressions into the first-order condition gives
\begin{align*}
     \nabla_\pi L(\pi,\lambda)\big|_{\pi=\pi_t}
    ={}&
    \alpha_t C
    +
    \log(K_t\oslash K_{t-1})
    \nonumber\\
    &+
    \left[
        \theta_t
        \log\left(
            \left(
                u_{t-1}\odot K_{t-1}v_{t-1}
            \right)
            \oslash p_t
        \right)
        +
        \log(u_t\oslash u_{t-1})
    \right]
    \one^\top
    \nonumber\\
    &+
    \one
    \left[
        \lambda
        +
        \log(v_t\oslash v_{t-1})
    \right]^\top.
\end{align*}

Hence, \( \nabla_\pi L(\pi,\lambda)\big|_{\pi=\pi_t}=0\) is ensured by the three conditions
\begin{equation*}
\begin{cases}
    \alpha_t C
    +
    \log(K_t\oslash K_{t-1})
    =0,
    \\[2mm]
    \displaystyle
    \theta_t
    \log\left(
        \left(
            u_{t-1}\odot K_{t-1}v_{t-1}
        \right)
        \oslash p_t
    \right)
    +
    \log(u_t\oslash u_{t-1})
    =0,
    \\[3mm]
    \displaystyle
    \lambda
    +
    \log(v_t\oslash v_{t-1})
    =0.
\end{cases}
\end{equation*}

The first and third equations immediately yield
\begin{equation*}
    K_t
    =
    K_{t-1}\odot\exp(-\alpha_t C),
\end{equation*}
and
\begin{equation*}
    \lambda
    =
    -\log(v_t\oslash v_{t-1}).
\end{equation*}

It remains to simplify the second equation. Substituting
\(
    p_t
    =
    p\odot
    u_{t-1}^{\frac{\theta_t-1}{\theta_t}},
\)
we obtain
\begin{align*}
    0
    ={}&
    \theta_t
    \log\left(
        \left(
            u_{t-1}\odot K_{t-1}v_{t-1}
        \right)
        \oslash
        \left(
            p\odot
            u_{t-1}^{\frac{\theta_t-1}{\theta_t}}
        \right)
    \right)
    +
    \log(u_t\oslash u_{t-1})
    \\
    ={}&
    \theta_t
    \log\left(
        u_{t-1}^{1/\theta_t}
        \odot
        \left(
            (K_{t-1}v_{t-1})\oslash p
        \right)
    \right)
    +
    \log(u_t\oslash u_{t-1})
    \\
    ={}&
    \log u_{t-1}
    +
    \theta_t
    \log\left(
        (K_{t-1}v_{t-1})\oslash p
    \right)
    +
    \log u_t
    -
    \log u_{t-1}.
\end{align*}
Therefore,
\[
    \log u_t
    =
    -\theta_t
    \log\left(
        (K_{t-1}v_{t-1})\oslash p
    \right),
\]
and hence
\begin{equation*}
    u_t
    =
    \left(
        p\oslash(K_{t-1}v_{t-1})
    \right)^{\theta_t}.
\end{equation*}

We have therefore shown that the OMD update is equivalent to
\begin{equation*}
\begin{cases}
    K_t
    =
    K_{t-1}\odot\exp(-\alpha_t C),
    \\[2mm]
    \displaystyle
    u_t
    =
    \left(
        p\oslash(K_{t-1}v_{t-1})
    \right)^{\theta_t},
    \\[3mm]
    \displaystyle
    \lambda
    =
    -\log(v_t\oslash v_{t-1}).
\end{cases}
\end{equation*}
Recall that the PAS column update
$v_t=q\oslash(K_t^\top u_t)$ is chosen to enforce
$\pi_t^\top\one=q$. Thus $\pi_t$ also satisfies the constraint
in the OMD problem, and the calculation above shows that
$\widetilde\pi_t=\pi_t$. Thus, this is precisely the PAS update and the potential-preserving iteration admits an OMD interpretation in which the factor \(\theta_t\) accounts explicitly for the temperature transfer.

\subsubsection{One-step KL identity}\label{sec:one-step-KL-id}

Let \(\pi\in\Gamma(p,q)\). 
As the second marginal is exact, for any \(\lambda\in \R^n\) we get \(\langle \lambda,\pi^\top\one-q\rangle=0\) and the first-order optimality condition associated to the dual problem from equation~\eqref{eq:OMG_Lagrangian} becomes
\begin{equation*}
\left\langle
\alpha_t C
+
\theta_t\log((\pi_{t-1}\one)\oslash p_t)\one^{\top}
+
\log(\pi_t\oslash\pi_{t-1}),
\,\pi-\pi_t
\right\rangle
= 0.
\end{equation*}
Using the definition \(p_t=p\odot u_{t-1}^{(\theta_t-1)/\theta_t}\) with \(\theta_t = \frac{\beta_t}{\beta_{t-1}}\), this becomes
\begin{align*}
\left\langle
\alpha_t C
+
\theta_t\log((\pi_{t-1}\one)\oslash p)\one^{\top}
+
\log(\pi_t\oslash\pi_{t-1}),
\,\pi-\pi_t
\right\rangle
=
(\theta_t-1)
\left\langle
\log(u_{t-1})\one^{\top},
\,\pi-\pi_t
\right\rangle .
\end{align*}

Moreover,
\[
\log(\pi_{t-1})
=
\log(u_{t-1})\one^{\top}
-\beta_{t-1}C
+\one\log(v_{t-1})^\top,
\]
giving
\begin{align*}
&\left\langle
\alpha_t C
+
\theta_t\log((\pi_{t-1}\one)\oslash p)\one^{\top}
+
\log(\pi_t\oslash\pi_{t-1}),
\,\pi-\pi_t
\right\rangle
\\
\nonumber
&\qquad=
(\theta_t-1)
\left\langle
\log(\pi_{t-1})
+\beta_{t-1}C
-\one\log(v_{t-1})^\top,
\,\pi-\pi_t
\right\rangle .
\end{align*}

Since \(\theta_t-1=\frac{\alpha_t}{\beta_{t-1}}\) and \(\pi\) and \(\pi_t\) have the same second marginal \(q\), we can simplify to obtain
\begin{align}
\label{eq:OMD_step_KL_equality}
\theta_t\underbrace{\left\langle
\log((\pi_{t-1}\one)\oslash p)\one^{\top},
\pi-\pi_t
\right\rangle}_{T_1}
+
\underbrace{\left\langle\log(\pi_t\oslash\pi_{t-1}),
\pi-\pi_t
\right\rangle}_{T_2}
=
\frac{\alpha_t}{\beta_{t-1}}\underbrace{
\left\langle
\log(\pi_{t-1}),
\,\pi-\pi_t
\right\rangle}_{T_3}.
\end{align}

As \(\pi\one=p\), using the three-point identity~\eqref{eq:three_point_identity} together with the definition \(\mathcal H(\pi)=\sum_{ij}\pi_{ij}\log\pi_{ij}\) gives
\[
\begin{aligned}
T_1 = \left\langle
\log (\pi_{t-1}\one \oslash p),
p-\pi_t\one
\right\rangle
={}&
\KL(\pi_t\one\mid\pi_{t-1}\one)
-\KL(\pi_t\one\mid p)
-\KL(p\mid\pi_{t-1}\one),
\\
T_2 = \left\langle
\log(\pi_t \oslash \pi_{t-1}),
\pi-\pi_t
\right\rangle
={}&
\KL(\pi\mid\pi_{t-1})
-\KL(\pi\mid\pi_t)
-\KL(\pi_t\mid\pi_{t-1}),
\\
\langle\log(\pi_t\oslash \one\one^\top),\pi-\pi_t\rangle
={}&
\KL(\pi\mid\one\one^\top)
-\KL(\pi\mid\pi_t)
-\KL(\pi_t\mid\one\one^\top)
\\
={}&
\mathcal H(\pi)
-\KL(\pi\mid\pi_t)
-\mathcal H(\pi_t).
\end{aligned}
\]
Hence, subtracting \(T_2\) from the last equality above we get
\[
T_3 = \left\langle
\log(\pi_{t-1}),
\,\pi-\pi_t
\right\rangle 
=
\mathcal H(\pi)
-\mathcal H(\pi_t)
-\KL(\pi\mid\pi_{t-1})
+\KL(\pi_t\mid\pi_{t-1}).
\]
Substituting the obtained identities for \(T_1,T_2,T_3\) into~\eqref{eq:OMD_step_KL_equality} gives
\begin{align*}
\theta_t&\left[\KL(\pi_t\one\mid\pi_{t-1}\one)
-\KL(\pi_t\one\mid p)
-\KL(p\mid\pi_{t-1}\one)
\right]
\\
&+
\KL(\pi\mid\pi_{t-1})
-\KL(\pi\mid\pi_t)
-\KL(\pi_t\mid\pi_{t-1})
\\=&
(\theta_t-1)
\left[\mathcal H(\pi)
-\mathcal H(\pi_t)
-\KL(\pi\mid\pi_{t-1})
+\KL(\pi_t\mid\pi_{t-1})\right].
\end{align*}
Cancelling the terms associated with \(\KL(\pi\mid\pi_{t-1})\) and \(\KL(\pi_t\mid\pi_{t-1})\) on both sides gives
\begin{align*}
\theta_t&\left[\KL(\pi_t\one\mid\pi_{t-1}\one)
-\KL(\pi_t\one\mid p)
-\KL(p\mid\pi_{t-1}\one)
\right]
-\KL(\pi\mid\pi_t)
\\=&
\left(\frac{\beta_t}{\beta_{t-1}}-\frac{\beta_t}{\beta_t}\right)
\left[\mathcal H(\pi)
-\mathcal H(\pi_t)\right]
-\theta_t\KL(\pi\mid\pi_{t-1})
+\theta_t\KL(\pi_t\mid\pi_{t-1}).
\end{align*}
Hence, dividing by \(\beta_t\) and rearranging, we obtain the one-step KL identity
\begin{align}
\label{eq:energy_balance}
\frac1{\beta_{t-1}}\Big[&
\KL(\pi_{t}\one\mid p)+\KL(p\mid \pi_{t-1}\one)
+\KL(\pi_t\mid\pi_{t-1})-\KL(\pi_{t}\one\mid \pi_{t-1}\one)
\Big]\\
\qquad
&=
\frac{1}{\beta_{t-1}}\KL(\pi\mid\pi_{t-1})-\frac{1}{\beta_t}\KL(\pi\mid\pi_t)
+\left(\frac1{\beta_{t-1}}-\frac1{\beta_t}\right)
\big[\mathcal H(\pi_t)-\mathcal H(\pi)\big].
\nonumber
\end{align}

\subsubsection{Subsequence convergence}
\label{subsubsec:subseq_conv}

We now show that there exists a subsequence $(t_j)$ such that $\pi_{t_j}\longrightarrow\bar \pi$, where \(\bar \pi\) is a solution of \(\OT(p,q)\).

As \(\bar \pi \in\Gamma(p,q)\), we can substitute it for \(\pi\) in the one-step KL identity~\eqref{eq:energy_balance}. Then, using the non-negativity of $\KL$~\eqref{eq:KL_non_neg} and the data processing inequality~\eqref{eq:data_processing_ineq}, we obtain
\begin{align}
\label{eq:kl_bar_pi_ineq}
\frac1{\beta_{t-1}}
\KL(\pi_{t}\one\mid p)
&=
-\frac1{\beta_{t-1}}\big[\overbrace{\KL(p\mid \pi_{t-1}\one)}^{\geq 0}
+\overbrace{\KL(\pi_t\mid\pi_{t-1})-\KL(\pi_{t}\one\mid \pi_{t-1}\one)}^{\geq 0}\big]\\
&+\frac{1}{\beta_{t-1}}\KL(\bar \pi\mid\pi_{t-1})-\frac{1}{\beta_t}\KL(\bar \pi\mid\pi_t)
+\left(\frac1{\beta_{t-1}}-\frac1{\beta_t}\right)
\big[\mathcal H(\pi_t)-\mathcal H(\bar \pi)\big].
\nonumber\\
&\leq \frac{1}{\beta_{t-1}}\KL(\bar \pi\mid\pi_{t-1})-\frac{1}{\beta_t}\KL(\bar \pi\mid\pi_t)
+\left(\frac1{\beta_{t-1}}-\frac1{\beta_t}\right)
\big[\mathcal H(\pi_t)-\mathcal H(\bar \pi)\big].
\nonumber
\end{align}
Since \(\pi_t\) and \(\bar \pi\) are probability measures on $mn$ entries,
\[
\mathcal H(\pi_t)-\mathcal H(\bar \pi)\le\log(mn).
\]
Summing from \(t=t_0+1\) to \(T\), with \(1\le t_0<T\), we obtain
\begin{align}
\label{eq:sub_sum_ineq}
\sum_{t=t_0+1}^{T}
\frac{1}{\beta_{t-1}}
\KL(\pi_t\one\mid p)
\leq
\frac{1}{\beta_{t_0}}
\KL(\bar\pi\mid\pi_{t_0})
+
\left(
\frac{1}{\beta_{t_0}}
-
\frac{1}{\beta_T}
\right)
\log(mn).
\end{align}

By assumption, the series
\(\sum_{t=t_0+1}^{\infty}\beta_{t-1}^{-1}\) diverges, whereas
\(
\sum_{t=t_0+1}^{\infty}
\beta_{t-1}^{-1}\KL(\pi_t\one\mid p)
\)
is finite by the estimate above. Consequently,
\(\liminf_t\KL(\pi_t\one\mid p)=0\), and there exists a sequence of indices
\((t_j)\) such that
\[
\KL(\pi_{t_j}\one\mid p)\longrightarrow0.
\]
Pinsker's inequality~\eqref{eq:pinskers_ineq} gives
\(
\|\pi_{t_j}\one-p\|_1^2
\leq2\KL(\pi_{t_j}\one\mid p)\to0.
\)
Since \(\pi_{t_j}^\top\one=q\), by compactness, a further subsequence, still denoted
$(\pi_{t_j})$, converges to a plan $\bar\pi$. Because
\(\beta_{t_j}\to\infty\), Lemma~\ref{lem:B1_Chizat} identifies its
limit as a solution of \(\OT(p,q)\).

\subsubsection{Convergence of the dual value}

Recall the dual values
\[
    L_t=p^\top f_t+q^\top g_t,
    \qquad
    f_t=\beta_t^{-1}\log u_t,
    \quad
    g_t=\beta_t^{-1}\log v_t.
\]
We prove that $(L_t)_{t\ge1}$ is nondecreasing and
converges to $\OT(p,q)$.

Using
\(
\log\pi_t
=\beta_t\bigl(f_t\one^\top+\one g_t^\top-C\bigr),
\)
and writing $\bar \pi\in\Gamma(p,q)$ for any fixed optimal solution of \(\OT(p,q)\), we have
\begin{align}
\label{eq:minor_dual_relation}
\frac{1}{\beta_t}\KL(\bar \pi\mid\pi_t)
=\langle C,\bar \pi\rangle-L_t+\frac{1}{\beta_t}\mathcal H(\bar \pi).
\end{align}
Substituting this identity into~\eqref{eq:kl_bar_pi_ineq}, we obtain
\begin{align*}
\frac1{\beta_{t-1}}
\KL(\pi_{t}\one\mid p)
\leq& \frac{1}{\beta_{t-1}}\KL(\bar \pi\mid\pi_{t-1})-\frac{1}{\beta_t}\KL(\bar \pi\mid\pi_t)
+\left(\frac1{\beta_{t-1}}-\frac1{\beta_t}\right)
\big[\mathcal H(\pi_t)-\mathcal H(\bar \pi)\big]
\\
=&
\langle C,\bar \pi\rangle - L_{t-1} + \frac{1}{\beta_{t-1}}\mathcal H(\bar \pi)
-\langle C,\bar \pi\rangle + L_t - \frac{1}{\beta_t}\mathcal H(\bar \pi) 
\\
&+\left(\frac1{\beta_{t-1}}-\frac1{\beta_t}\right)
\big[\mathcal H(\pi_t)-\mathcal H(\bar \pi)\big]
\\
=&L_t - L_{t-1} +\left(\frac1{\beta_{t-1}}-\frac1{\beta_t}\right)
[\mathcal H(\pi_t)].
\end{align*}
Using the non-negativity of \(\KL\)~\eqref{eq:KL_non_neg},
and
\(\mathcal H(\pi_t)\leq0\), as all entries of \(\pi_t\) are non-negative and at most 1,
we obtain
\begin{align*}
L_t
\geq
L_{t-1} + \frac1{\beta_{t-1}}
\KL(\pi_{t}\one\mid p)
-\left(\frac1{\beta_{t-1}}-\frac1{\beta_t}\right)\mathcal H(\pi_t)
\geq L_{t-1}.
\end{align*}
Since \(\pi_t^\top\one=q\), the entries of \(\pi_t\) are at most one. 
As
\(f_{t,i}+g_{t,j}-C_{ij}=\beta_t^{-1}\log\pi_{t,ij}\leq0\),
then
\[
    L_t
    =\sum_{i}f_{t,i}p_i+\sum_{j}g_{t,j}q_j
    =\sum_{i,j}(f_{t,i}+g_{t,j})\bar \pi_{ij}
    \le\langle C,\bar \pi\rangle
    = \OT(p,q).
\]
The previous section provides a subsequence such that
\[
\pi_{t_j}\longrightarrow\bar\pi,
\qquad
\bar\pi\in\operatorname*{argmin}_{\pi\in\Gamma(p,q)}
\langle C,\pi\rangle,
\]
and by equation~\eqref{eq:minor_dual_relation}, we have
\[
\OT(p,q)-L_{t_j}
=\frac{1}{\beta_{t_j}}\KL(\bar\pi\mid\pi_{t_j})
-\frac{1}{\beta_{t_j}}\mathcal H(\bar\pi),
\]
Since \(\KL(\bar\pi\mid\pi_{t_j})\to0\) and \(\beta_{t_j}\to \infty\), we have \(L_{t_j}\to\OT(p,q)\). Hence, by monotonicity of the sequence \((L_t)\), we must in fact have
\[
L_t\uparrow\OT(p,q).
\]

\subsubsection{Identification of the limiting plan}
\label{proof:Thm2_uniqueness}

Recall that $\pi^\star$ is the unique minimiser of
$\mathcal H$ on $\mathcal S_{\rm OT}$, defined in
\eqref{C1}. We show that every accumulation point
$\bar\pi\in\mathcal S_{\rm OT}$ equals $\pi^\star$.

Consider a subsequence \(\pi_{t_j} \to \bar \pi\) for some \(\bar \pi \in \mathcal S_{\rm OT}\). Then \(\langle C,\bar \pi\rangle = \langle C,\pi^\star\rangle\) and, by the previous section, \(
L_t\uparrow\OT(p,q).
\) Hence, by \eqref{eq:minor_dual_relation} we have
\[
\KL(\pi^\star\mid\pi_{t_j}) - \mathcal H(\pi^\star) = \KL(\bar \pi\mid\pi_{t_j}) - \mathcal H(\bar \pi).
\]
By nonnegativity of \(\KL\), and since \(\KL(\bar \pi \mid \pi_{t_j})\to 0\), taking limits gives
\[
\mathcal H(\pi^\star) \geq \mathcal H(\bar \pi).
\]
By minimality of \(\pi^\star\), we also have
\[
\mathcal H(\pi^\star) \leq \mathcal H(\bar \pi).
\]
Strict convexity of \(\mathcal H\) gives
\[
\pi^\star = \bar \pi.
\]
Thus, every accumulation point in \(\mathcal S_{\rm OT}\) equals \(\pi^\star\), and the compactness of \(\Gamma(*,q)\) gives
\[
\KL(\pi^\star\mid\pi_{t_j}) \to 0.
\]

\subsection{Proof of Proposition~\ref{prop:negative}}
\label{app:proof_prop}

Let
\[
 S=\sum_{t\geq0}\beta_t^{-1}<\infty,\qquad
 C=\begin{bmatrix}0&1\\1&0\end{bmatrix},\qquad
 q=[1/2,1/2],\qquad
 p=[1/2+\eta,1/2-\eta],
\]
where $\eta>0$ is small enough that
$S\log\!\bigl((1/2+\eta)/(1/2-\eta)\bigr)<1$, which holds, for example, whenever $\eta<(4S+2)^{-1}$.
It suffices to show that the PAS iterates converge to
$\operatorname{diag}(1/2,1/2)$, whose source marginal is $[1/2,1/2]\neq p$.
Initialise with zero potentials and let
\[
 \ell=\log\frac{p_1}{p_2},\qquad
 x_t=\frac1{\beta_t}\log\frac{u_{t,1}}{u_{t,2}},\qquad
 y_t=\frac1{\beta_t}\log\frac{v_{t,1}}{v_{t,2}}.
\]
The first row update gives $x_1=\ell/\beta_0$.  Because the column projection enforces the symmetric target $q=[1/2,1/2]$ exactly, its log-ratio satisfies
$-x_t\leq y_t\leq0$.  Substituting this bound in the next row projection gives the exact one-step control
\[
 0<x_{t+1}\leq x_t+\frac{\ell}{\beta_t}.
\]
Consequently,
\[
 0<x_t\leq \ell\sum_{s=0}^{t-1}\frac1{\beta_s}
 \leq \ell S<1.
\]
The within-column ratios of the post-projection plan are therefore
\[
 \frac{\pi_{t,11}}{\pi_{t,21}}
   =e^{\beta_t(1+x_t)}\longrightarrow\infty,
 \qquad
 \frac{\pi_{t,12}}{\pi_{t,22}}
   =e^{-\beta_t(1-x_t)}\longrightarrow0.
\]
Since both column sums equal $1/2$, it follows that
$\pi_t\to\operatorname{diag}(1/2,1/2)$.  The limit violates the source constraint whenever $\eta>0$.

\subsection{Calculation for Example~\ref{ex:non_conv_example}}
\label{sec:example_calculation}
Let $p=(a,1-a)^\top$, where
$a\in(0,1)\setminus\{1/2\}$, and take $q=(1)$ and
$C=(0,1)^\top$. Set $\beta_0=1$ and
\[
\beta_t=2^k
\qquad\text{for }2^k\le t<2^{k+1},\quad k\ge0.
\]

For \(\theta>0\), write
\[
Z_\theta=a^\theta+(1-a)^\theta.
\]
Since there is only one target point, the PAS update gives
\[
u_{t,i}
=p_i^{\theta_t}e^{\beta_t C_i}v_{t-1}^{-\theta_t},
\qquad
v_t=\frac{v_{t-1}^{\theta_t}}{Z_{\theta_t}}.
\]
Consequently,
\[
(\pi_t\one)_i=\frac{p_i^{\theta_t}}{Z_{\theta_t}},
\qquad
f_t=C+\frac{\log p}{\beta_{t-1}}
       -g_{t-1}\mathbf1,
\qquad
g_t=g_{t-1}-\frac{\log Z_{\theta_t}}{\beta_t}.
\]
The prescribed schedule gives $\theta_{2^k}=2$ for $k\ge1$
and $\theta_t=1$ at all other indices $t\ge1$.
Since $Z_1=1$ and $Z_2=a^2+(1-a)^2$, we obtain
\[
    \pi_{2^k}\one
    =\frac{1}{Z_2}\bigl(a^2,(1-a)^2\bigr)^\top,
    \qquad
    \pi_{2^k+1}\one=p,
    \qquad k\ge1.
\]
These two constant subsequences have distinct limits
because $a\ne1/2$. Since every first marginal equals
one of these two vectors, they are exactly the two
accumulation points of $(\pi_t\one)$.
The expressions for $f_t$ and $g_t$ above give
\(
f_{t,i}+g_t
=C_i+\frac{\log((\pi_t\one)_i)}{\beta_t}.
\)
Hence,
\[
L_t
=1-a+\frac{1}{\beta_t}\sum_{i=1}^2p_i\log((\pi_t\one)_i),
\]
Since $(\pi_t\one)_i\in\{p_i,p_i^2/Z_2\}\subset(0,1)$,
there exists $M=M(a)>0$ such that
\[
    0\le (1-a)-L_t\le\frac{M}{\beta_t}\longrightarrow0.
\]
Together with the monotonicity established in
Theorem~\ref{thm:positive}, this gives
$L_t\uparrow\OT(p,q)=1-a$.

\subsection{Proof of Theorem~\ref{thm:full_convergence}}
\label{app:thm_full_conv_proof}

We will first derive a one-step marginal KL bound that separates the usual Sinkhorn decrease from the contribution of annealing. Using the local geometry of KL, we will then absorb this contribution into the negative KL term up to a summable error. Together with Theorem~\ref{thm:positive}, the resulting estimate
will imply $\pi_t\one\to p$. Lemma~\ref{lem:B1_Chizat} and the identification of the limiting plan in Section~\ref{proof:Thm2_uniqueness} will then imply
convergence of the entire sequence to $\pi^\star$.

Throughout the proof we will use the notation defined in~\eqref{eq:notation}, namely
\(\alpha_t= \beta_t-\beta_{t-1}\) and \(\theta_t=\frac{\beta_t}{\beta_{t-1}}\). 
We will also make use of the KL projection $\mathcal{P}_p$ onto the set of matrices with first marginal $p$. Writing $\widehat u=p\oslash(K_{t-1}v_{t-1})$ for the row update
at inverse temperature $\beta_{t-1}$, the resulting intermediate
plan is
\begin{equation}\label{eq:row_proj}
    \mathcal{P}_p(\pi_{t-1})
    :=\operatorname{diag}(\widehat u)K_{t-1}
      \operatorname{diag}(v_{t-1})
    =\bigl[(p\oslash(\pi_{t-1}\one))\one^\top\bigr]
      \odot\pi_{t-1}.
\end{equation}
It satisfies $\mathcal{P}_p(\pi_{t-1})\one=p$ and hence has unit mass.

\subsubsection{One-step marginal inequality}
A careful re-examination of steps taken in Section~\ref{sec:one-step-KL-id} to arrive at  \eqref{eq:OMD_step_KL_equality}, which is obtained from the OMD optimality condition, reveals that it in fact holds for every comparison plan $\pi \in \Gamma(*,q)$.
In particular, substituting $\pi=\pi_{t-1}$ results in the following equality
\begin{align*}
\left\langle
\theta_t\log((\pi_{t-1}\one)\oslash p)\one^{\top}
+
\log(\pi_t\oslash\pi_{t-1}),
\,\pi_{t-1}-\pi_t
\right\rangle
=
\frac{\alpha_t}{\beta_{t-1}}
\left\langle
\log(\pi_{t-1}),
\,\pi_{t-1}-\pi_t
\right\rangle.
\end{align*}
We now isolate the contribution from the change in inverse temperature. By the definition of the row projection in~\eqref{eq:row_proj}, we have
\[
    \log\mathcal{P}_p(\pi_{t-1})
    =
    \log\pi_{t-1}
    -
    \log((\pi_{t-1}\one)\oslash p)\one^\top.
\]
Writing $\theta_t=1+\frac{\alpha_t}{\beta_{t-1}}$ in the preceding
equality and moving the part of the first inner product
proportional to $\frac{\alpha_t}{\beta_{t-1}}$ to the right-hand side
therefore gives
\begin{align*}
&\left\langle
    \log((\pi_{t-1}\one)\oslash p),
    \pi_{t-1}\one-\pi_t\one
\right\rangle
+
\left\langle
    \log(\pi_t\oslash\pi_{t-1}),
    \pi_{t-1}-\pi_t
\right\rangle
\\
&\qquad=
\frac{\alpha_t}{\beta_{t-1}}
\left\langle
    \log\mathcal{P}_p(\pi_{t-1}),
    \pi_{t-1}-\pi_t
\right\rangle.
\end{align*}
For the first inner product, the three-point identity
\eqref{eq:three_point_identity}, applied to the row marginals,
gives
\[
\begin{aligned}
\left\langle
    \log((\pi_{t-1}\one)\oslash p),
    \pi_{t-1}\one-\pi_t\one
\right\rangle
={}&
\KL(\pi_{t-1}\one\mid p)
-\KL(\pi_t\one\mid p)
\\
&+\KL(\pi_t\one\mid\pi_{t-1}\one).
\end{aligned}
\]
This contains the difference between the successive marginal
KL discrepancies. For the second inner product, the same
identity yields
\[
\begin{aligned}
\left\langle
    \log(\pi_t\oslash\pi_{t-1}),
    \pi_{t-1}-\pi_t
\right\rangle
&=
-\left\langle
    \log(\pi_t\oslash\pi_{t-1}),
    \pi_t-\pi_{t-1}
\right\rangle
\\
&=
-\KL(\pi_t\mid\pi_{t-1})
-\KL(\pi_{t-1}\mid\pi_t).
\end{aligned}
\]

Substituting these identities and rearranging for
$\KL(\pi_t\one\mid p)$, we obtain
\begin{align*}
\KL(\pi_t\one\mid p)
={}&
\KL(\pi_{t-1}\one\mid p)
-\KL(\pi_{t-1}\mid\pi_t)
\\
&-
\bigl[
    \KL(\pi_t\mid\pi_{t-1})
    -\KL(\pi_t\one\mid\pi_{t-1}\one)
\bigr]
\\
&+
\frac{\alpha_t}{\beta_{t-1}}
\left\langle
    \log\mathcal{P}_p(\pi_{t-1}),
    \pi_t-\pi_{t-1}
\right\rangle.
\end{align*}
The bracketed term is nonnegative by the data processing
inequality~\eqref{eq:data_processing_ineq}. Hence, for $t\ge2$,
\begin{equation}
\label{eq:initial_KL_step_bound}
\KL(\pi_t\one\mid p)
\le
\KL(\pi_{t-1}\one\mid p)
-\KL(\pi_{t-1}\mid\pi_t)
+\frac{\alpha_t}{\beta_{t-1}}\mathcal A_t,
\end{equation}
where
\[
\mathcal A_t
:=\left\langle
\log\mathcal{P}_p(\pi_{t-1}),\pi_t-\pi_{t-1}
\right\rangle.
\]
The first two terms on the right-hand side give the usual
marginal KL decrease for standard Sinkhorn, while the term
including $\mathcal A_t$ arises from annealing.
Without this additional term, $\KL(\pi_t\one\mid p)$ would
be nonincreasing, and Theorem~\ref{thm:positive} would imply
its convergence to zero.
We will therefore now show that the annealing contribution
can be absorbed into the negative KL term up to a summable
error.

\subsubsection{Controlling the annealing term}

In this section, we will show that
\[
|\mathcal A_t|
\le \sqrt{2B\KL(\pi_{t-1}\mid\pi_t)},
\]
for some constant $B>0$ independent of $t$.
In what follows, for notational brevity, we will often
suppress dependence on $t$ when defining auxiliary terms.

The crucial first step is to normalise the columns of
$\mathcal{P}_p(\pi_{t-1})$ without changing $\mathcal A_t$.
To this end, set
\begin{equation}\label{eq:normalisation}
    \ell_j:=\sum_i\mathcal{P}_p(\pi_{t-1})_{ij},
    \qquad
    z_{ij}:=\frac{\mathcal{P}_p(\pi_{t-1})_{ij}}{\ell_j}.
\end{equation}
Since $\pi_t$ and $\pi_{t-1}$ have the same column marginal $q$,
we have
\[
\begin{aligned}
\mathcal A_t
&=\langle\log z,\pi_t-\pi_{t-1}\rangle
  +\sum_j\log\ell_j
    \underbrace{\sum_i(\pi_{t,ij}-\pi_{t-1,ij})}_{=0}\\
&=\langle\log z,\pi_t-\pi_{t-1}\rangle.
\end{aligned}
\]
Next, we use the Cauchy--Schwarz inequality to obtain
\begin{align*}
\left|\mathcal A_t\right|
&=
\left|
\sum_{ij}
(\sqrt{\pi_{t}}-\sqrt{\pi_{t-1}})_{ij}
(\sqrt{\pi_{t}}+\sqrt{\pi_{t-1}})_{ij}
\log z_{ij}
\right|\\
&\leq
\left(
\sum_{ij}(\sqrt{\pi_{t}}-\sqrt{\pi_{t-1}})_{ij}^2
\right)^{1/2}
\left(
\sum_{ij}
(\sqrt{\pi_{t}}+\sqrt{\pi_{t-1}})_{ij}^2
(\log z_{ij})^2
\right)^{1/2}.
\end{align*}
The bound involves two factors: the Hellinger distance \citep{tsybakov_introduction_2009} and the square root of a weighted second moment of $\log(z)$. This inequality reflects the quadratic local geometry of KL and is reminiscent of those derived in \citep{bolley_weighted_2005} and used in \citep{ghosal_convergence_2026} to derive bounds for standard Sinkhorn. Since the square of the Hellinger distance is bounded by KL \citep{tsybakov_introduction_2009}, for the first term, we obtain
\[
\left(
\sum_{ij}(\sqrt{\pi_{t}}-\sqrt{\pi_{t-1}})_{ij}^2
\right)^{1/2}
\le
\sqrt{\KL(\pi_{t-1}\mid \pi_{t})}.
\]
It thus remains to bound the second factor by a constant independent
of $t$. To achieve that, we first use Young's inequality on \((\sqrt{\pi_{t}}+\sqrt{\pi_{t-1}})_{ij}^2\) to split it into two parts
\begin{align}
    \label{eq:two_sums_to_bound}
    \sum_{ij}
    (\sqrt{\pi_{t}}+\sqrt{\pi_{t-1}})_{ij}^2
    (\log z_{ij})^2
    \leq
    2
    \sum_{ij}
    \pi_{t,ij}
    (\log z_{ij})^2
    +
    2
    \sum_{ij}
    \pi_{t-1,ij}
    (\log z_{ij})^2.
\end{align}
We first bound $\sum_{ij}\pi_{t-1,ij}(\log z_{ij})^2$.
The definitions of $\mathcal{P}_p$ in~\eqref{eq:row_proj}
and $z_{ij}$ in~\eqref{eq:normalisation} give
\[
    \ell_j z_{ij}
    =\mathcal{P}_p(\pi_{t-1})_{ij}
    =\frac{p_i\pi_{t-1,ij}}{(\pi_{t-1}\one)_i}.
\]
Rearranging, we obtain
\[
\pi_{t-1,ij} = \frac{z_{ij}\ell_j \sum_k \pi_{t-1, ik}}{p_i} \leq \frac{z_{ij}\ell_j}{p_{\min}},
\]
as \(\sum_k \pi_{t-1,ik}\le 1\) and \(p_{\min}:=\min_i p_i>0\). Now, as 
\(
x(\log x)^2\leq\frac4{e^2}
\) 
for all \(x\in (0,1]\)
and
\(\sum_j\ell_j=1\), 
we obtain
\begin{align*}
\sum_{ij}\pi_{t-1,ij}(\log z_{ij})^2
\le
\frac1{p_{\min}}
\sum_j\ell_j\sum_i z_{ij}(\log z_{ij})^2
\le
\frac1{p_{\min}}
\sum_j\ell_j\frac{4m}{e^2}
=
\frac{4m}{p_{\min}e^2}
\eqcolon
C_1.
\end{align*}
This bounds the second sum in~\eqref{eq:two_sums_to_bound}.
For the first sum, we express the PAS update in terms of
$\mathcal{P}_p(\pi_{t-1})$. The update consists of raising
the entries of this intermediate plan to the power $\theta_t$,
followed by a column projection onto $q$. Thus
\begin{align}
\label{eq:first_sum_cancelation}
\pi_{t,ij}
&=
q_j
\frac{\bigl(\mathcal{P}_p(\pi_{t-1})_{ij}\bigr)^{\theta_t}}
{\sum_k\bigl(\mathcal{P}_p(\pi_{t-1})_{kj}\bigr)^{\theta_t}}
=
q_j\frac{\ell_j^{\theta_t}z_{ij}^{\theta_t}}
{\ell_j^{\theta_t}\sum_k z_{kj}^{\theta_t}}
=
q_j\frac{z_{ij}^{\theta_t}}{Z_j},
\end{align}
where $Z_j:=\sum_k z_{kj}^{\theta_t}$.

As \(z_{ij}^{\theta_t}\le z_{ij}\), in order to use the same inequality 
\(
x(\log x)^2\leq\frac4{e^2}
\) 
that we used for the second sum, we need to bound \(Z_j\). From the assumption \eqref{AS2} we know
\(
\sum_{t=1}^{\infty}
\left(\frac{\beta_t - \beta_{t-1}}{\beta_{t-1}}\right)^2<\infty,
\)
we thus immediately get that
\(\theta_t \to 1\) and \(1\leq \theta_t \leq \Theta\) for some fixed \(\Theta > 1\).

As \(Z_{j} = \sum_{k=1}^m z_{k j}^{\theta_t}\), using Jensen's inequality
we have 
\[
\frac1m\sum_{k=1}^m z_{kj}^{\theta_t}
\ge
\left(\frac1m\sum_{k=1}^m z_{kj}\right)^{\theta_t}
=
\left(\frac1m\right)^{\theta_t},
\]
giving
\[
Z_j
\ge m^{1-\theta_t}
\geq
m^{1-\Theta}
=
\frac{1}{C_z},
\]
for \(C_z=m^{\Theta-1}\). As
\[
Z_j\ge\frac1{C_z}
\qquad
\text{and}
\qquad
z_{ij}^{\theta_t}\le z_{ij},
\]
we have
\(
\pi_{t,ij}\le C_zq_jz_{ij}.
\)
This gives the bound
\[
\begin{aligned}
\sum_{ij}\pi_{t,ij}(\log z_{ij})^2
\le
C_z\sum_jq_j
\sum_i z_{ij}(\log z_{ij})^2
\le
C_z\sum_jq_j\frac{4m}{e^2}
=\frac{4mC_z}{e^2}
\eqcolon
C_2.
\end{aligned}
\]
This shows why the normalisation in~\eqref{eq:normalisation}
was needed, as without it the denominator
in~\eqref{eq:first_sum_cancelation} would only admit a bound
\(\sum_k(\mathcal{P}_p(\pi_{t-1})_{kj})^{\theta_t}
\geq m^{1-\Theta}\ell_j^{\theta_t}\).
The bound on the total mass alone does not prevent individual column masses \(\ell_j\) from being very small, which would make us unable to obtain a constant independent of \(\pi_{t-1}\).
Putting everything together, with \(B=2C_1+2C_2\) we obtain the final bound
\begin{align*}
\left|\mathcal{A}_t\right|
&\leq
\left(
\sum_{ij}(\sqrt{\pi_{t}}-\sqrt{\pi_{t-1}})_{ij}^2
\right)^{1/2}
\left(
\sum_{ij}
(\sqrt{\pi_{t}}+\sqrt{\pi_{t-1}})_{ij}^2
(\log z_{ij})^2
\right)^{1/2}\\
&\leq
\sqrt{2B\KL(\pi_{t-1}\mid \pi_{t})}.
\end{align*}

\subsubsection{Summable errors and full convergence}

The preceding bound on $|\mathcal A_t|$ and Young's inequality
$ab\le\frac12(a^2+b^2)$, applied to
${a=\sqrt{\KL(\pi_{t-1}\mid\pi_t)}}$ and
${b=\sqrt{2B}\,\alpha_t/\beta_{t-1}}$, give
\begin{align*}
\frac{\alpha_t}{\beta_{t-1}}\mathcal A_t
&\le
\frac{\alpha_t}{\beta_{t-1}}
\sqrt{2B\KL(\pi_{t-1}\mid\pi_t)}\\
&\le
\frac12\KL(\pi_{t-1}\mid\pi_t)
+B\left(\frac{\alpha_t}{\beta_{t-1}}\right)^2.
\end{align*}
Substituting into~\eqref{eq:initial_KL_step_bound}
and using $\beta_t/\beta_{t-1}\le\Theta$, we obtain
\begin{equation}
\label{ineq:one_step_KL}
\KL(\pi_t\one\mid p)
\le
\KL(\pi_{t-1}\one\mid p)
-\frac12\KL(\pi_{t-1}\mid\pi_t)
+C_0\left(\frac{\alpha_t}{\beta_t}\right)^2,
\end{equation}
where $C_0:=B\Theta^2$.
Since $\alpha_t=\beta_t-\beta_{t-1}$, this is the
claimed one-step estimate.

Since $\beta_t\ge\beta_{t-1}$, assumption~\eqref{AS2} implies
\[
\sum_{t=1}^{\infty}\left(\frac{\alpha_t}{\beta_t}\right)^2
< \infty,
\]
so we can define the tail-corrected sequence
\[
D_t:=\KL(\pi_t\one\mid p)
+C_0\sum_{j=t+1}^{\infty}\left(\frac{\alpha_j}{\beta_j}\right)^2.
\]
By inequality~\eqref{ineq:one_step_KL},
\[
\begin{aligned}
D_{t+1}
&=\KL(\pi_{t+1}\one\mid p)
+C_0\sum_{j=t+2}^{\infty}\left(\frac{\alpha_j}{\beta_j}\right)^2\\
&\leq\KL(\pi_t\one\mid p)
+C_0\left(\frac{\alpha_{t+1}}{\beta_{t+1}}\right)^2
+C_0\sum_{j=t+2}^{\infty}\left(\frac{\alpha_j}{\beta_j}\right)^2\\
&=D_t.
\end{aligned}
\]
Thus, \((D_t)\) is a nonnegative, nonincreasing sequence and therefore converges. Since
\[
\sum_{j=k}^{\infty}\left(\frac{\alpha_j}{\beta_j}\right)^2\longrightarrow0
\qquad\text{as }k\to\infty,
\]
\(\KL(\pi_t\one\mid p)\) also converges.
By Theorem~\ref{thm:positive}, we have
$\liminf_t\KL(\pi_t\one\mid p)=0$.
Since the sequence converges, it follows that
\[
\KL(\pi_t\one\mid p)\longrightarrow0.
\]

The plans $\pi_t$ lie in the compact probability simplex and always satisfy $\pi_t^\top\one=q$. By Pinsker's inequality
\(
\pi_t\one\to p\text{ in }\ell^1
\). Hence, every accumulation point also has source marginal $p$. Lemma~\ref{lem:B1_Chizat} identifies every such accumulation point as a solution of the corresponding OT problem. 
The argument in Section~\ref{proof:Thm2_uniqueness} identifies every such accumulation point with $\pi^\star$, the maximum-entropy optimal transport plan. Compactness therefore gives $\pi_t\to\pi^\star$, and hence $\KL(\pi^\star\mid\pi_t)\to0$.

We note that the marginal convergence argument above is a
deterministic instance of the Robbins--Siegmund convergence
theorem \citep{robbins_convergence_1971} 
with the corresponding scalar convergence lemmas stated in
\cite[Lemma~5.31]{bauschke_convex_2017} and
\cite[Lemma~3.6]{franci_convergence_2022}. The argument uses the same
summable-error principle as quasi-Fej\'er analysis
\cite{combettes_quasi_2001}.
Indeed, \eqref{ineq:one_step_KL} has the form
\[
d_{t+1}\le d_t-cS_t+e_t,
\]
with
\[
\begin{aligned}
d_t&:=\KL(\pi_t\one\mid p),
& S_t&:=\KL(\pi_t\mid\pi_{t+1}),\\
e_t&:=C_0\left(\frac{\alpha_{t+1}}{\beta_{t+1}}\right)^2,
& c&:=\frac12.
\end{aligned}
\]
These quantities are nonnegative, and \eqref{AS2} implies
$\sum_t e_t<\infty$. The sequence used above is precisely
$D_t=d_t+\sum_{j=t}^{\infty}e_j$.

\subsection{Proof of Theorem~\ref{thm:convergence_rate}}

We will first obtain an $O(T^{-1})$ bound on the minimum marginal
KL discrepancy over $\lfloor T/2\rfloor<t\le T$.
We will then use the one-step estimate from
Theorem~\ref{thm:full_convergence} to control its subsequent
increase and obtain the same rate at iteration $T$.
Finally, Pinsker's inequality and the rounding estimate yield
the marginal $\ell^1$ and transport cost bounds.

\subsubsection{A bound on the minimum marginal KL}
Let \(T\geq2\) and \(a=\lfloor T/2\rfloor\), so that \(a\ge1\).
We use the telescoping estimate~\eqref{eq:sub_sum_ineq} of Section~\ref{subsubsec:subseq_conv}, starting at \(t=a\). Then
\begin{align}
\sum_{t=a+1}^{T}
\frac{1}{\beta_{t-1}}
\KL(\pi_t\one\mid p)
\leq
\frac{1}{\beta_{a}}
\KL(\pi^\star\mid\pi_{a})
+
\left(
\frac{1}{\beta_{a}}
-
\frac{1}{\beta_T}
\right)
\log(mn),
\end{align}
where \(\pi^\star\) is the unique maximum-entropy plan defined in~\eqref{C1}.
By Theorem~\ref{thm:full_convergence}, we have 
\(
\KL(\pi^\star\mid\pi_t)\to0,
\)
which implies that 
\(
\sup_{t\ge1} \KL(\pi^\star\mid\pi_t)<A_0,
\)
for some fixed \(A_0 > 0\) and
\[
\frac{1}{\beta_{T-1}}
\sum_{t=a+1}^{T}
\KL(\pi_t\one\mid p)
\leq
\sum_{t=a+1}^{T}
\frac{1}{\beta_{t-1}}
\KL(\pi_t\one\mid p)
\leq
\frac{A_0 + \log(mn)}{\beta_a}.
\] 
Since $\beta_t = (t+1)^{\kappa}$ for some \(0 < \kappa \leq 1\) we have
\[
\frac{\beta_{T-1}}{\beta_a} \leq \left(\frac{2T}{T+1}\right)^\kappa \leq 2
\]
Hence we obtain
\[
\sum_{t=a+1}^{T}
\KL(\pi_t\one\mid p)
\leq
2(A_0 + \log(mn)).
\] 
The sum contains \(T-a\geq T/2\) terms, so
\[
\min_{\lfloor T/2\rfloor<t\le T}
\KL(\pi_t\one\mid p)
\le \frac{A_1}{T},
\]
where \(A_1 = 4(A_0 + \log(mn))\).
\subsubsection{A marginal KL bound at the final iterate}
Choose
\(
s\in\{\lfloor T/2\rfloor+1,\dots,T\}
\)
such that
\(
\KL(\pi_s\one\mid p)
=
\min_{\lfloor T/2\rfloor<t\le T}\KL(\pi_t\one\mid p).
\)
Then, by inequality~\eqref{ineq:one_step_KL},
\[
\KL(\pi_T\one\mid p)
\le
\KL(\pi_s\one\mid p)+
C_0\sum_{t=s+1}^T\left(\frac{\alpha_t}{\beta_t}\right)^2.
\]
As \(\beta_t=(t+1)^\kappa\) and \(x^\kappa\) is concave for \(0\leq\kappa\leq1\),
\[
\alpha_t = \beta_t - \beta_{t-1} = (t+1)^\kappa-t^\kappa\leq \kappa t^{\kappa-1}.
\]
Therefore
\[
\begin{aligned}
\left(\frac{\beta_t-\beta_{t-1}}{\beta_t}\right)^2
\leq \frac{\kappa^2t^{2\kappa-2}}{(t+1)^{2\kappa}}
\leq \frac{\kappa^2}{t^2}.
\end{aligned}
\]
Because \(s>T/2\),
\[
\sum_{t=s+1}^T\left(\frac{\alpha_t}{\beta_t}\right)^2
\le
\kappa^2\sum_{t=s+1}^T\frac1{t^2}
\le
\kappa^2 \frac{T}{2} \left(\frac{2}{T}\right)^2
\leq
\frac{2}{T},
\]
as \(\kappa \leq 1\).
Combining these bounds, we obtain
\begin{align}
\KL(\pi_T\one\mid p)
\le
\KL(\pi_s\one\mid p)+
C_0\sum_{j=s+1}^T\left(\frac{\alpha_j}{\beta_j}\right)^2
\leq  
\frac{A_1}{T}
+
\frac{2C_0}{T}
=
\frac{A}{T}
\end{align}
for \(A = A_1 + 2C_0\) and \(T\ge2\). Enlarging \(A\), if necessary, so that \(A\ge\KL(\pi_1\one\mid p)\) also covers \(T=1\) and proves the first claim.
\subsubsection{Marginal $\ell^1$ and transport cost bounds}
Pinsker's inequality then gives
\[
\|\pi_T\one-p\|_1\le\sqrt{2A}\,T^{-1/2}.
\]
For the rounding map $R$ defined in~\eqref{eq:rounding_the_plan}, the cost estimate in \citep[Lemma~2.3]{chizat_annealed_2024} gives
\begin{align*}
\langle C,R(\pi_T)\rangle-\OT(p,q) 
&\leq \frac{1}{\beta_T}\log(mn) + 4\left(\|\pi_T\one-p\|_1+ \|\pi_T^\top\one-q\|_1\right)\|C\|_{\mathrm{osc}}
\\
&\leq \log(mn)T^{-\kappa} + \sqrt{32A}\|C\|_{\mathrm{osc}}T^{-1/2}
\end{align*}
where \(\|C\|_{\mathrm{osc}}=\max_{ij}C_{ij}-\min_{ij}C_{ij}\) is the cost oscillation and \(\|\pi_T^\top\one-q\|_1 = 0\) as the second marginal of any \(\pi_t\) is always exactly \(q\). For \(B = \log(mn) + \sqrt{32A}\|C\|_{\mathrm{osc}}\) depending only on $p,q,C$ and $\kappa$, we get the final bound
\begin{align*}
\langle C,R(\pi_T)\rangle-\OT(p,q) 
\leq B(T^{-\kappa} + T^{-1/2}).
\end{align*}

\bibliographystyle{plainnat}
\bibliography{refs}

\clearpage
\appendix

\section{Annealing schedules and state transfer in practice}
\label{app:libraries}

In this appendix, we examine two practical aspects of annealed Sinkhorn methods: how the state transfer between successive temperatures is handled in leading computational optimal transport libraries, and the choice of piecewise-constant annealing schedules.

\subsection{State transfer in computational OT libraries}

As discussed in Section~\ref{sec:PAS}, the scaling and potential parametrisations are equivalent at fixed temperature, but retaining scalings or potentials across a temperature change gives different annealed iterations. Since both rules occur in widely used computational OT libraries, we record the annealed iteration realised by each implementation. We hope that this may provide a useful perspective on how annealed Sinkhorn is represented in computational implementations.

Recall that we call annealed Sinkhorn \emph{scaling-preserving} (SAS) when $(u,v)$ are retained as the temperature changes, and \emph{potential-preserving} (PAS) when the corresponding dual potentials
\[
    (f,g)=\bigl(\beta^{-1}\log u,\beta^{-1}\log v\bigr)
\]
are retained instead. The descriptions below refer to repository snapshots accessed on 23 November 2026, and the corresponding commit hashes are reported for reproducibility.

\paragraph{POT.}
At commit
\href{https://github.com/PythonOT/POT/tree/367355d757223e39989e384b7a93d0537b1547ea}
{\texttt{367355d}}
of POT \citep{flamary_pot_2021}, the routine
\href{https://github.com/PythonOT/POT/blob/367355d757223e39989e384b7a93d0537b1547ea/ot/bregman/_sinkhorn.py\#L1515-L1557}
{\texttt{sinkhorn\_epsilon\_scaling}}
uses the stabilised Sinkhorn solver at each value of $\varepsilon$. Between successive values of $\varepsilon$, it retains the absorbed dual variables in cost units, i.e. the physical dual potentials, and therefore realises PAS.

\paragraph{GeomLoss.}
At commit
\href{https://github.com/jeanfeydy/geomloss/tree/00e493f36bd6cc8471526e59afc07193a1926e47}
{\texttt{00e493f}}
of GeomLoss \citep{feydy_interpolating_2019}, the
\href{https://github.com/jeanfeydy/geomloss/blob/00e493f36bd6cc8471526e59afc07193a1926e47/src/geomloss/ot/_abstract_solvers/sinkhorn_ot.py\#L244-L272}
{$\varepsilon$-scaling Sinkhorn loop}
is formulated in dual potentials in cost units and retains these as $\varepsilon$ changes, hence realising PAS. GeomLoss uses symmetrised Sinkhorn updates, together with optional debiasing and multiscale operations, so the full iteration differs from the alternating PAS iteration analysed here.

\paragraph{OTT-JAX.}
At commit
\href{https://github.com/ott-jax/ott/tree/b8fbb2b34367a7cc9ddf19f6496435cb19ef0932}
{\texttt{b8fbb2b}}
of OTT-JAX \citep{cuturi_optimal_2022}, both state-transfer rules occur. With an
\href{https://github.com/ott-jax/ott/blob/b8fbb2b34367a7cc9ddf19f6496435cb19ef0932/src/ott/geometry/epsilon_scheduler.py}
{$\varepsilon$ schedule}
and
\href{https://github.com/ott-jax/ott/blob/b8fbb2b34367a7cc9ddf19f6496435cb19ef0932/src/ott/solvers/linear/sinkhorn.py\#L803-L897}
{\texttt{lse\_mode=True}},
the iterated variables are the dual potentials $(f,g)$, which are retained across temperature changes; this gives PAS. With \texttt{lse\_mode=False}, termed \emph{kernel mode}, the iterated variables are instead the multiplicative scalings $(u,v)$, giving SAS. The documentation presents these modes as two numerical implementations of the same fixed-temperature Sinkhorn iteration: kernel mode uses matrix multiplication and is faster per iteration, whereas log-sum-exp mode is more stable for small $\varepsilon$.
Under annealing, however, retaining the variables of the chosen representation fixes the state-transfer rule, and hence the annealed iteration.
The numerical parametrisation and the state-transfer rule are independent choices.
For example, PAS may be implemented entirely in scaling variables by applying, at a change $\beta\to\beta'$,
\[
    u^+=\widehat u^{\,\beta'/\beta},
    \qquad
    v^+=\widehat v^{\,\beta'/\beta},
\]
while SAS may equally be expressed in potential variables through
\[
    f^+=\frac{\beta}{\beta'}\,\widehat f,
    \qquad
    g^+=\frac{\beta}{\beta'}\,\widehat g.
\]

\paragraph{OptimalTransport.jl.}
At commit
\href{https://github.com/JuliaOptimalTransport/OptimalTransport.jl/tree/7c0a19624ee24a416f375c81ca7aa4d8fb622de7}
{\texttt{7c0a196}}
of OptimalTransport.jl \citep{juliaoptimaltransport_optimaltransportjl_2026},
\href{https://github.com/JuliaOptimalTransport/OptimalTransport.jl/blob/7c0a19624ee24a416f375c81ca7aa4d8fb622de7/src/entropic/sinkhorn_epsscaling.jl\#L34-L103}
{\texttt{SinkhornEpsilonScaling}}
applies $\varepsilon$-scaling to several Sinkhorn variants. 
The
\href{https://github.com/JuliaOptimalTransport/OptimalTransport.jl/blob/7c0a19624ee24a416f375c81ca7aa4d8fb622de7/src/entropic/sinkhorn_gibbs.jl\#L14-L43}
{\texttt{SinkhornGibbs}}
formulation retains the multiplicative scalings $(u,v)$ between temperature levels and hence realises SAS.

The
\href{https://github.com/JuliaOptimalTransport/OptimalTransport.jl/blob/7c0a19624ee24a416f375c81ca7aa4d8fb622de7/src/entropic/sinkhorn_stabilized.jl\#L103-L132}
{\texttt{SinkhornStabilized}}
formulation uses a mixed parametrisation. Writing the represented potentials as
\[
    f=\alpha+\varepsilon\log u,
    \qquad
    g=\gamma+\varepsilon\log v,
\]
it retains both the absorbed components $(\alpha,\gamma)$ in cost units and the residual scalings $(u,v)$. Since these components transform differently when $\varepsilon$ changes, the resulting rule is in general neither SAS nor PAS and depends on the preceding absorption steps.

\paragraph{FlashSinkhorn.}
At commit
\href{https://github.com/ot-triton-lab/flash-sinkhorn/tree/888540b4f2675e9646c10a49f5c042136b6d513f}
{\texttt{888540b}}
of FlashSinkhorn \citep{ye_flashsinkhorn_2026}, the
\href{https://github.com/ot-triton-lab/flash-sinkhorn/blob/888540b4f2675e9646c10a49f5c042136b6d513f/torch-ext/flash_sinkhorn/sinkhorn_solvers.py\#L647-L718}
{$\varepsilon$-scaling solver}
uses symmetric Sinkhorn updates in shifted dual potentials. These potentials
are retained as $\varepsilon$ changes, so the state-transfer rule is PAS.
As in GeomLoss, the symmetric within-temperature iteration differs from the
alternating PAS iteration analysed here.

\subsection{Geometric piecewise-constant schedules}
A standard $\varepsilon$-scaling scheme uses geometrically decreasing temperatures and several Sinkhorn iterations at each temperature. The libraries above use closely related schedules: OTT-JAX and OptimalTransport.jl anneal geometrically, GeomLoss does so in $\varepsilon^{1/p}$ and POT decreases $\varepsilon-\varepsilon_{\star}$ geometrically towards a prescribed terminal temperature. Corollary~\ref{cor:piecewise_constant} gives a simple worst-case condition under which such stagewise schemes retain the subsequential convergence guarantee of Theorem~\ref{thm:positive}.
\begin{corollary}[Piecewise-constant schedules]
\label{cor:piecewise_constant}
Let $(\bar\beta_\ell)_{\ell\ge1}$ be a positive, nondecreasing
sequence of inverse temperatures tending to infinity, and let
$M_\ell\in\mathbb N_{>0}$. Define $(\beta_t)_{t\ge0}$ by
successively repeating each $\bar\beta_\ell$ exactly $M_\ell$
times. If
\[
\sum_{\ell=1}^{\infty}M_\ell\bar\beta_\ell^{-1}=\infty,
\]
then the PAS iterates of Algorithm \ref{alg:PAS} admit a subsequence
converging to the maximum-entropy optimal transport plan $\pi^\star$.
\end{corollary}
\begin{proof}
The schedule $(\beta_t)$ is positive and nondecreasing, tends
to infinity, and satisfies
\[
\sum_{t=0}^{\infty}\beta_t^{-1}
=\sum_{\ell=1}^{\infty}M_\ell\bar\beta_\ell^{-1}
=\infty.
\]
The conclusion follows from Theorem~\ref{thm:positive}.
\end{proof}

For fixed geometric growth $\bar\beta_{\ell+1}/\bar\beta_\ell=\rho>1$, the relative increment at every temperature transition equals $\rho-1$, so \eqref{AS2} fails regardless of the stage lengths. Theorem~\ref{thm:full_convergence} therefore does not strengthen this corollary to full convergence. Thus geometric cooling remains admissible for the subsequential guarantee provided the number of Sinkhorn iterations per temperature grows sufficiently rapidly. We next quantify the corresponding iteration counts and attainable regularisation levels.

Consider a geometric inverse-temperature schedule
\[
    \bar\beta_\ell=\beta_0\rho^\ell,
    \qquad
    \bar\varepsilon_\ell=\bar\beta_\ell^{-1},
    \qquad
    \rho>1,
\]
and let $M_\ell$ be the number of fixed-temperature Sinkhorn iterations performed at the $\ell$th temperature level. Corollary~\ref{cor:piecewise_constant}
then requires
\[
    \sum_{\ell\geq1}\frac{M_\ell}{\bar\beta_\ell}
    =
    \frac1{\beta_0}
    \sum_{\ell\geq1}M_\ell\rho^{-\ell}
    =
    \infty.
\]
Hence $M_\ell$ must compensate for the geometric growth of $\bar \beta_{\ell}$.
A simple sufficient choice is
\[
    M_\ell=\left\lceil\frac{\rho^\ell}{\ell}\right\rceil,
\]
for which
\[
    \frac{M_\ell}{\bar\beta_\ell}
    \geq \frac{1}{\beta_0\ell}.
\]
Table~\ref{tab:geometric-inner-iterations} gives the terminal temperature, the final-level iteration count $M_L$ and the cumulative count $\sum_{\ell=1}^L M_{\ell}$ for this schedule.
We take $\beta_0=10$ and $\max C-\min C=1$.

\begin{table}[h]
\centering
\small
\setlength{\tabcolsep}{7pt}
\renewcommand{\arraystretch}{1.12}
\begin{tabular}{@{}rrrrr@{}}
\toprule
$\rho$
& $L$
& $\bar\varepsilon_L$
& $M_L$
& $\sum_{\ell=1}^L M_\ell$ \\
\midrule
$1.25$
& $31$
& $9.90\times10^{-5}$
& $33$
& $212$ \\

$2$
& $10$
& $9.77\times10^{-5}$
& $103$
& $240$ \\

$2$
& $12$
& $2.44\times10^{-5}$
& $342$
& $769$ \\

$2$
& $15$
& $3.05\times10^{-6}$
& $2{,}185$
& $4{,}756$ \\

$2$
& $18$
& $3.81\times10^{-7}$
& $14{,}564$
& $31{,}127$ \\

$4$
& $5$
& $9.77\times10^{-5}$
& $205$
& $303$ \\

$4$
& $7$
& $6.10\times10^{-6}$
& $2{,}341$
& $3{,}327$ \\

$4$
& $8$
& $1.53\times10^{-6}$
& $8{,}192$
& $11{,}519$ \\
\bottomrule
\end{tabular}
\caption{Iteration counts for the geometric schedule
$\bar\beta_\ell=\beta_0\rho^\ell$ with
$M_\ell=\lceil\rho^\ell/\ell\rceil$ and $\beta_0=10$.
The final column gives $\sum_{\ell=1}^L M_\ell$.}
\label{tab:geometric-inner-iterations}
\end{table}

For example, with $\rho=2$, a terminal temperature of order $10^{-4}$ is
reached using $240$ Sinkhorn iterations in total, while temperatures of order
$10^{-6}$ and $10^{-7}$ require approximately $5\times10^3$ and
$3\times10^4$ iterations, respectively.
Thus the sufficient worst-case condition remains computationally moderate over a substantial finite range, although $M_{\ell}$ eventually inherits the geometric growth of $\bar \beta_{\ell}$.

\end{document}